\documentclass[11pt,twoside,a4paper]{article}

\usepackage[cp1250]{inputenc} 
\usepackage{lmodern}
\usepackage{amsmath}
\usepackage{indentfirst}
\usepackage{amsmath,amssymb}
\usepackage{microtype}
\usepackage{graphicx}
\usepackage{amsthm}
\DisableLigatures{encoding = *, family = * }
\usepackage{fancyhdr}
\usepackage{pstricks,graphicx}
\usepackage{amssymb}
\usepackage{float}
\usepackage{hyperref}
\usepackage{bm}
\usepackage{tabularx}
\usepackage{authblk}
\usepackage[title]{appendix}

\usepackage{color}
\usepackage{geometry}

\newgeometry{tmargin=2.5cm, bmargin=2.5cm, headheight=14.5pt, inner=3cm, outer=2.5cm} 

\definecolor{color02}{rgb}{0.00,0.00,1.00}
\newcommand{\TT}{\mathsf{T}}

\newcommand{\mrm}[1]{\mathrm{#1}}

\renewcommand{\vec}[1]{\bm{#1}}
\newcommand{\mat}[1]{\mathbf{#1}}

\newcommand{\x}[1]{\vec{x}_{#1}}

\newcommand{\E}{\mathbb{E}}

\newcommand{\cN}{\mathcal{N}}

\usepackage{amsthm}
\theoremstyle{theorem}
\newtheorem{thm}{Theorem}
\newtheorem{lem}{Lemma}

\newcommand{\eqnref}[1]{Eq.~(\ref{#1})}
\newcommand{\eqnsref}[2]{Eqs.~(\ref{#1}-\ref{#2})}

\begin{document}
\title{An Approximate Bayesian Approach to Optimal Input Signal Design for System Identification}

\definecolor{lime}{HTML}{A6CE39}
\DeclareRobustCommand{\orcidicon}{
    \begin{tikzpicture}
    \draw[lime, fill=lime] (0,0) 
    circle [radius=0.16] 
    node[white] {{\fontfamily{qag}\selectfont \tiny ID}};
    \draw[white, fill=white] (-0.0625,0.095) 
    circle [radius=0.007];
    \end{tikzpicture}
    \hspace{-2mm}
}
\foreach \x in {A, ..., Z}{\expandafter\xdef\csname orcid\x\endcsname{\noexpand\href{https://orcid.org/\csname orcidauthor\x\endcsname}
            {\noexpand\orcidicon}}
}

\author[1, 2]{\small Piotr Bania\thanks{pba@agh.edu.pl}}
%\email{pba@agh.edu.pl}
\author[1]{\small Anna W\'{o}jcik\thanks{annawoj@student.agh.edu.pl}}

\affil[1]{\footnotesize AGH University of Krak\'{o}w, Faculty of Electrical Engineering, Automation, Computer Science and Biomedical Engineering, Department of Automation and Robotics, al.~A.~Mickiewicza 30, 30-059 Krak\'{o}w, Poland.}
\affil[2]{\footnotesize Institute of Physics, Polish Academy of Sciences, Aleja Lotnik\'{o}w 32/46, 02-668 Warsaw, Poland.}
%\affil[3]{\footnotesize Department of Astrophysical Sciences, Princeton University, 4 Ivy Lane, Princeton, NJ 08544-1001 USA}
%\affil[4]{\footnotesize Department of Astronomy, Yale University, New Haven, CT 06520, USA}
%\affil[5]{\footnotesize Einstein and Spitzer Fellow}

\maketitle
%\date{\today}
\begin{abstract}
 The design of informatively rich input signals is essential for accurate system identification, yet classical Fisher-information-based methods are inherently local and often inadequate in the presence of significant model uncertainty and nonlinearity. This paper develops a Bayesian approach that uses the mutual information (MI) between observations and parameters as the utility function. To address the computational intractability of the MI, we maximize a tractable MI lower bound. The method is then applied to the design of an input signals for the identification of quasi-linear stochastic dynamical systems. Evaluating the MI lower bound requires inversion of large covariance matrices whose dimensions scale with the number of data points $N$. To overcome this problem, an algorithm that reduces the dimension of the matrices to be inverted by a factor of $N$ is developed, making the approach feasible for long experiments. The proposed Bayesian method is compared with the average D-optimal design method, a semi-Bayesian approach, and its advantages are demonstrated. The effectiveness of the proposed method is further illustrated through four examples, including atomic sensor models, where the input signals that generates large MI are especially important for reducing the estimation error.
\end{abstract}
%\begin{keywords} a,b,c
%\end{keywords}


\section{Introduction}
The design of informative input signals is the cornerstone of modern system identification. Without properly chosen excitation, even advanced estimation algorithms may fail to provide accurate parameter estimates, leading to unreliable prediction and control. Classical references such as \cite{GoodwinPayne1977}, \cite{Ljung1999}, \cite{soderstrom1989} and reviews by \cite{Pronzato2008}, \cite{HuanJagalurMarzouk2024}, \cite{RainforthEtAl2024}, emphasize that identification is not only a matter of statistical estimation but also of experimental design, where the input signal determines the achievable information content. Optimal experimental design (OED) methods therefore play a crucial role in practical applications. Traditionally, OED has relied on the Fisher Information Matrix (FIM), with criteria such as D or A optimality widely used due to their computational efficiency and asymptotic guarantees \cite{Pronzato2008}, \cite{fedorov1997}. However, FIM-based approaches are inherently local, relying on linearization and asymptotic normality. They may thus be fragile in scenarios with large model uncertainty or strongly non-linear stochastic dynamics.

A natural alternative is the Bayesian approach, which evaluates an experiment through its \emph{expected information gain}, typically quantified by the mutual information between model parameters and observations \cite{HuanJagalurMarzouk2024}, \cite{RainforthEtAl2024}, \cite{Lindley1956}, \cite{ArKi1971}, \cite{Chaloner1995}, \cite{Ryan2015}. Bayesian design has several advantages: it is globally valid over the parameter space, it naturally incorporates prior information, and it is applicable to non-linear, stochastic models. Its main drawback is computational intractability, since mutual information requires high-dimensional integration over parameters and observations. Since the general case is challenging and difficult to solve, in this article we focus on models that can be represented in the form
\begin{equation}
\vec{Y}=\vec{F}(\vec{\theta}, \vec{U})+\vec{Z},
\label{eq:model_0}
\end{equation}
where $\vec{\theta}\in\lbrace \vec{\theta}_1, \ldots, \vec{\theta}_r\rbrace$ is a parameter with prior distribution $P(\vec{\theta}=\vec{\theta}_j)=p_{0,j}$. The noise $\vec{Z}$ is conditionally normal, that is, $p(\vec{Z}|\vec{\theta})=\cN(\vec{Z},0,\mat{S}(\vec{\theta}, \vec{U}))$, and $\vec{U}$ is a design variable. For this class of models, the density of the observations $\vec{Y}$ is a finite Gaussian mixture of the form $p(\vec{Y})=\sum\limits_{j=1}^r p_{0,j}\cN(\vec{Y},\vec{F}(\vec{\theta}_j, \vec{U}), \mat{S}(\vec{\theta}_j, \vec{U}))$. Within these Gaussian mixtures, the mutual information between $\vec{Y}$ and $\vec{\theta}$ can be estimated from below by the effective and tractable pairwise-distance-based lower bound given by Kolchinsky and Tracey \cite{Kolchinsky2017}, \cite{Kolchinsky_arxiv}. We maximize this bound to achieve an approximate optimal design parameter $\vec{U}$ and then generalize the method to a parameter space of continuum cardinality. In particular, we show how to treat the Gaussian prior $p_0(\vec{\theta})=\cN(\vec{\theta},\vec{m}_{\theta}, \mat{S}_{\theta})$ and prior distributions with compact support.

In this work, we focus on \emph{quasi-linear systems}, namely stochastic dynamical systems that are linear in the state variables but non-linear in the control variables. Such systems occur ubiquitously in science and engineering. In quantum mechanics, Hamiltonians and Lindblad dissipators depend on external control fields such as laser intensities, magnetic fields, or gate voltages, leading to a non-linear dependence on the control \cite{AltafiniTicozzi2012}, \cite{DongPetersen2010}. In chemical processes, flow rates directly determine reaction speeds in continuous stirred tank reactors \cite{Friedly1972}, \cite{Lorenz2007Discrimination}. In thermal plants, convection coefficients scale non-linearly with flow, giving rise to quasi-linear heat transfer dynamics \cite{Friedly1972}. This broad applicability makes quasi-linear systems a natural and important class for advanced input design. However, there is a notable lack of Bayesian design methods and software tools tailored to this class of systems. Therefore, in this paper, we address this gap by developing such a method. Specifically, we show that a finite sequence of observations generated by a quasi-linear system can always be expressed in the form of the model \eqref{eq:model_0}, and we provide an effective algorithm for calculating the lower bound on mutual information.

The study reported here constitutes a substantial and far-reaching extension of the initial results presented in \cite{Bania2019}, as well as related research reported in \cite{BanBar2016}, \cite{Ban2018}, \cite{BarBanAl2017}. The article's main contributions can be summarized as follows. We first introduce an Information-Theoretic Lower Bound (ITB) on the estimation error of any estimator, \cite{CoverThomas2006}, \cite{Lee2022}, and briefly discuss its relation to the Bayesian Cram\'{e}r-Rao Bound (BCRB), \cite{Trees1968}, \cite{Efroimovich1979}, \cite{VanTreesBell2007_BayesianBounds} \cite{Lee2022}.
We conclude that maximizing mutual information is superior to maximizing Bayesian or classical Fisher information, which is consistent with the arguments presented in \cite{Lindley1956}, \cite{Chaloner1995}, \cite{Ryan2015}, \cite{HuanJagalurMarzouk2024}. Building on this result, we introduce a novel Bayesian design method for model \eqref{eq:model_0}. To address the intractability and computational complexity of direct mutual information evaluation, we discretize the parameter space and maximize the Kolchinsky-Tracey lower bound \cite{Kolchinsky2017}, \cite{Kolchinsky_arxiv}, and subsequently extend this approach to a parameter space of continuum cardinality. We then focus on the application to linear and quasi-linear system identification. Since the information-theoretic bound requires inversion of large covariance matrices of the observations, whose dimensions grow linearly with the number of data points $N$ (with $N \sim 10^3$-$10^6$ in applications), direct inversion is computationally problematic. To overcome this challenge, we develop an algorithm that reduces the dimension of the matrices needed for inversion by a factor of $N$, thus making the approach feasible for long experiments. The proposed Bayesian method is compared with the average D-optimal design \cite{GoodwinPayne1977}, \cite{Pronzato2008}, \cite{Jakowluk2024}, a semi-Bayesian method, and its advantages are demonstrated. The effectiveness of the method is further illustrated by four examples. The first two, intentionally elementary, highlights the effectiveness of our approach. The third and fourth examples are drawn from atomic sensor models, a domain where optimal input design is particularly critical. We analyze a controlled harmonic oscillator with stochastic disturbances as a paradigmatic atomic sensor model \cite{JimenezMatinez2018}, and a complex magnetometer model with a non-linear dependence of the system matrices on the input \cite{Truullinou2021}. In the latter case, we provide a simplified model, derive the optimal input, and demonstrate that it significantly outperforms the harmonic signal, which might otherwise be presumed optimal. Moreover, we show that the estimation error of the MAP estimator achieves the theoretical lower bound.

The article is organized as follows. Section~2 formulates the problem. Section~3 develops the approximate Bayesian solution for finite and infinite parameter spaces. Section~4 applies the method to quasi-linear systems. Section~5 compares the approach with classical design methods. Section~6 presents examples. Section~7 provides discussion and conclusions.

%**********************
\section{Formulation of the problem} Let us consider a family of models
\begin{equation}
\vec{Y}=\vec{F} (\vec{\theta},\vec{U})+\vec{Z},
\label{eq:model_01}
\end{equation}
where $\vec{Y}, \vec{Z}\in \mathbb{R}^{n_Y}$, $\vec{U}\in \mathbb{R}^{n_U}$, $\vec{\theta} \in \vec{\Theta} \subset \mathbb{R}^{n_{\theta}}$. The set $\vec{\Theta}$ will be called parameter space. Parameter $\vec{\theta}$ is unknown.  Prior distribution of $\vec{\theta}$ is denoted by $p_0$. Random variable $\vec{Z}$ is conditionally normal i.e. $p(\vec{Z}|\vec{\theta})=\cN(\vec{Z}, 0, \mat{S}(\vec{\theta}, \vec{U}))$, where $\mat{S}(\vec{\theta}, \vec{U})\in\mat{S}^{+}(n_Y)$, for all $\vec{\theta}\in\vec{\Theta}$, $\vec{U}\in\mathbb{R}^{n_U}$. Functions $\vec{F}$ and $\mat{S}$ are smooth. The variable $\vec{U}$ is called the design parameter or, in the context of dynamical systems, the input signal. The set of admissible signals is given by
\begin{equation}
\mathcal{U}_{ad}=\lbrace \vec{U}\in \mathbb{R}^{n_U}; |\vec{U}-\tilde{\vec{U}}|\leqslant\varrho \rbrace,
\label{eq:U_ad}
\end{equation}
where $\tilde{\vec{U}}$ is the given vector and $\varrho$ is maximal norm of the signal. We will also consider an alternative, and useful in some applications, definition of $\mathcal{U}_{ad}$:
\begin{equation}
\mathcal{U}_{ad}=\lbrace \vec{U}\in \mathbb{R}^{n_U}; \vec{U}_{min}\leqslant \vec{U} \leqslant \vec{U}_{max} \rbrace,
\label{eq:U_ad_1}
\end{equation}
where $\vec{U}_{min}, \vec{U}_{max}\in \mathbb{R}^{n_U}$ are fixed vectors. Under these assumptions and after applying the Bayes rule we get the likelihood, evidence and posterior distribution of $\vec{\theta}$:
\begin{equation}
 p(\vec{Y}| \vec{\theta},\vec{U})=\cN(\vec{Y},\vec{F}(\vec{\theta}, \vec{U}), \mat{S}(\vec{\theta},\vec{U})),
 \label{eq:likelihood_pdf}
\end{equation}
\begin{equation}
 p(\vec{Y}|\vec{U})=\int p_0(\vec{\theta})\cN(\vec{Y},\vec{F}(\vec{\theta}, \vec{U}), \mat{S}(\vec{\theta},\vec{U}))d\vec{\theta},
 \label{eq:pdfs_00}
\end{equation}
\begin{equation}
    p(\vec{\theta}|\vec{Y}, \vec{U})=\frac{p_0(\vec{\theta})p(\vec{Y}|\vec{\theta}, \vec{U})}{p(\vec{Y}|\vec{U})}.
\label{eq:pdfs_01}
\end{equation}
The Minimum Mean Squared Error (MMSE) estimator of $\vec{\theta}$ is then given by
\begin{equation}
    \hat{\vec{\theta}}(\vec{Y}, \vec{U})=\int \vec{\theta}p(\vec{\theta}|\vec{Y}, \vec{U})d\vec{\theta}.
    \label{eq:mse_estimator}
\end{equation}
To avoid the difficulties involved in calculating the integral \eqref{eq:mse_estimator}, instead of the MSE, the Maximum a Posteriori (MAP) estimator is typically used. Taking the negative logarithm of both sides of \eqref{eq:pdfs_01} and omitting the terms independent on $\vec{\theta}$, we get the following:
\begin{equation}
\mathcal{L}(\vec{\theta},\vec{Y},\vec{U})=\tfrac{1}{2}|\vec{Y}-\vec{F}(\vec{\theta}, \vec{U})|_{\mat{S}^{-1}(\vec{\theta}, \vec{U})}^2+\tfrac{1}{2}\ln|\vec{S}(\vec{\theta}, \vec{U})|-\ln p_0(\vec{\theta}).
\label{eq:log_likelihood_0}
\end{equation}
Thus, the MAP estimator of $\vec{\theta}$ is given by
\begin{equation}
\hat{\vec{\theta}}(\vec{Y},\vec{U})=\mathrm{arg}\min_{\vec{\theta}\in\vec{\Theta}}\mathcal{L}(\vec{\theta},\vec{Y},\vec{U}).
\label{eq:map_estimator}
\end{equation}
Estimators \eqref{eq:mse_estimator} or \eqref{eq:map_estimator} may be biased, so the Cram\'{e}r-Rao Bound (CRB) cannot be applied directly to them. However, a Bayesian version of the CRB exists and can be used to estimate the error of any, also biased, estimator \cite{Trees1968}, \cite{Efroimovich1979}, \cite{VanTreesBell2007_BayesianBounds}, \cite{Lee2022}. Let us introduce Bayesian information (BI):
\begin{equation}
    \mat{J}_B=\E_{p(\vec{\theta},\vec{Y}|\vec{U})}\left[\nabla_{\vec{\theta}}\mathcal{L}(\vec{\theta},\vec{Y},\vec{U})(\nabla_{\vec{\theta}}\mathcal{L}(\vec{\theta},\vec{Y},\vec{U})^{\TT})\right]=\mat{J}_D+\mat{J}_P,
    \label{eq:bayesian_inf_matrix}
\end{equation}
where
\begin{equation}
    \mat{J}_P=\E_{p_0(\vec{\theta})}\left[\nabla_{\vec{\theta}}\ln p_0(\vec{\theta})(\nabla_{\vec{\theta}}\ln p_0(\vec{\theta}))^{\TT}\right],
    \label{eq:J_P}
\end{equation}
is the fraction of BI associated to a prior and
\begin{equation}
    \mat{J}_D(\vec{U})=\E_{p(\vec{\theta}, \vec{Y}|\vec{U})}\left[\nabla_{\vec{\theta}}\ln p(\vec{Y}|\vec{\theta}, \vec{U})(\nabla_{\vec{\theta}}\ln p(\vec{Y}|\vec{\theta}, \vec{U}))^{\TT}\right],
    \label{eq:J_D}
\end{equation}
is part of BI provided by observations $\vec{Y}$. The matrix $\mat{J}_D$ is the Bayesian equivalent of the Fisher information matrix. The Fisher information matrix can be recovered from \eqref{eq:J_D} assuming $p_0(\vec{\theta})=\delta(\vec{\theta}-\vec{\theta}^{0})$, where $\vec{\theta}^{0}$ is the true value of the parameter. Let $\hat{\vec{\theta}}(\vec{Y},\vec{U})$ be any estimator of $\vec{\theta}$. Assuming a sufficiently regular prior, it can be proven that
\begin{equation}
    \E|\vec{\theta} - \hat{\vec{\theta}}(\vec{Y},\vec{U})|^2\geqslant\frac{n_{\theta}}{|\mat{J}_P+\mat{J}_D(\mat{U})|^{1/n_{\theta}}},
    \label{eq:van_trees}
\end{equation}
which is usually known as the Van Trees inequality or Bayesian Cram\'{e}r-Rao Bound (BCRB) \cite{Trees1968}, \cite{Bobrov1987}, \cite{Lee2022}[inequality (2.9), p. 17]. Formulas \eqref{eq:bayesian_inf_matrix}-\eqref{eq:J_D} are well defined only under rather restrictive assumptions. In particular, the joint density $p(\vec{Y},\vec{\theta})$ must be differentiable with respect to $\vec{\theta}$, and it must satisfy the regularity condition $\int \nabla_{\vec{\theta}} p(\vec{Y},\vec{\theta}) d\vec{Y} = 0$.
Moreover, the prior and likelihood distributions must guarantee the existence of the expectations in \eqref{eq:J_P} and \eqref{eq:J_D}. This excludes the uniform and many other useful prior distributions and considerably limits
the applicability of inequality \eqref{eq:van_trees} (see \cite{Lee2022} for details). Beyond inequality \eqref{eq:van_trees}, there exists a large class of Bayesian bounds reported in \cite{VanTreesBell2007_BayesianBounds}. Many of these bounds can serve as a utility function.

Probably one of the best design criterion is the Ziv-Zakai lower bound \cite{JeoDytso}. However, to compute this estimate, an additional, and rather complex, optimization sub-problem must be solved as shown in \cite{JeoDytso}. Therefore, due to the high computational complexity of the multivariate Ziv-Zakai bound, we will not consider it here.
Given the application-oriented focus of this article and following the arguments presented in \cite{HuanJagalurMarzouk2024}, \cite{RainforthEtAl2024}, we conclude that the entropy-based lower bound \cite{CoverThomas2006}[p. 255], \cite{Lee2022}[Sect. 2.2, p.16-17], is a reasonable optimality criterion and provides slightly tighter estimates than the BCRB (cf. \cite{JeoDytso}[Sect. V.D.]). To proceed, let us define the entropies of $\vec{Y}$ and $\vec{\theta}$ and the corresponding conditional entropies:
\begin{align}
\label{eq:entropies_00}
&H_{\vec{Y}}(\vec{U})=\E(-\ln p(\vec{Y}|\vec{U})),\\
&H_{\vec{\theta}}=\E(-\ln p_0(\vec{\theta})),\\
&H_{\vec{Y}|\vec{\theta}}(\vec{U})=\E(-\ln p(\vec{Y}|\vec{\theta}, \vec{U}))=\tfrac{1}{2}\int p_0(\vec{\theta})\ln \left((2\pi e)^{n_Y}|\mat{S}(\vec{\theta}, \vec{U})|\right)d\vec{\theta},\\
&H_{\vec{\theta}|\vec{Y}}(\vec{U})=\E(-\ln p(\vec{\theta}|\vec{Y}, \vec{U})).
\label{eq:entropies_01}
\end{align}
The mutual information (MI) between $\vec{\theta}$ and $\vec{Y}$ is defined as
\begin{equation}
I_{\vec{\theta};\vec{Y}}(\vec{U})=H_{\vec{\theta}}-H_{\vec{\theta}|\vec{Y}}(\vec{U})=H_{\vec{Y}}(\vec{U})-H_{\vec{Y}|\vec{\theta}}(\vec{U}).
\label{eq:MI}
\end{equation}
The following theorem establishes the ultimate limit of the estimation error expressed in terms of mutual information, demonstrating that the Bayesian Cram\'{e}r-Rao Bound (BCRB) does not constitute a fundamental limit.
\begin{thm}
\label{thm:theorem_1}
Let $\hat{\vec{\theta}}(\vec{Y},\vec{U})$ be any estimator of $\vec{\theta}$. 
Then the following inequalities hold: 
\begin{equation} 
\E|\vec{\theta} - \hat{\vec{\theta}}(\vec{Y},\vec{U})|^2\geqslant n_{\theta}(2\pi e)^{-1} e^{2n_{\theta}^{-1}(H_{\vec{\theta}}-I_{\vec{\theta};\vec{Y}}(\vec{U}))}\geqslant\frac{n_{\theta}}{|\mat{J}_P+\mat{J}_D(\vec{U})|^{1/n_{\theta}}}.
\label{eq:error_lb}
\end{equation}
\end{thm}
The proof is given in Appendix \ref{app:appendix_A}. The first of the inequalities \eqref{eq:error_lb} will be called the Information-Theoretic Lower Bound (ITB). The last part of \eqref{eq:error_lb} is known as the Efroimovich inequality \cite{Efroimovich1979}, \cite{Lee2022}[inequality (2.7), Ch. 2.2, p. 16].

To determine the optimal signal, one may maximize the MI, the determinant of $\mat{J}_B$ or the determinant of FIM. The latter corresponds to the classical design methods outlined in Sect.~\ref{sec:classical_methods}. However, the right-hand side of \eqref{eq:van_trees} generally underestimates the estimation error, and large values of $\mat{J}_B$ do not necessarily guarantee a small error. We illustrate this problem in Appendix \ref{app:appendix_B}.
In contrast, the ITB \eqref{eq:error_lb}, which plays a central role in our subsequent analysis, shows that maximizing $I_{\theta;Y}(\vec{U})$ is essential for reducing the estimation error and provides a more fundamental criterion than maximizing either the Bayesian or classical Fisher information. In particular, in the context of optimal experimental design and input signal design for system identification, maximizing the mutual information between the parameters and the observations constitutes the most principled optimality criterion, as it directly quantifies the amount of knowledge gained about the parameters \cite{HuanJagalurMarzouk2024}. Accordingly, we define the optimal signal as the solution of the following optimization problem:
\begin{equation}
\vec{U}^{\ast}=\mathrm{arg}\max_{\vec{U}\in\mathcal{U}_{ad}} I_{\vec{\theta};\vec{Y}}(\vec{U}).
\label{eq:optimization_problem}
\end{equation}
As $I_{\vec{\theta};\vec{Y}}$ is continuous and $\mathcal{U}_{ad}$ is compact, then \eqref{eq:optimization_problem} is well defined. After solving the task, the MMSE, MAP, or any other estimator can be used to determine $\vec{\theta}$.

Computing mutual information (MI) or its lower-bound remains a significant challenge. The analyses presented in the literature \cite{Chaloner1995}, \cite{Ryan2015}, \cite{HuanJagalurMarzouk2024}, \cite{RainforthEtAl2024}, show that this can be done in three main ways: ~1) using Monte Carlo or nested Monte Carlo (MC) simulations, \cite{HuanJagalurMarzouk2024}[Sect. 3.1.], \cite{RainforthEtAl2024}; ~2) applying variational lower-bound (VLB) estimates of the MI, \cite{HuanJagalurMarzouk2024}[Sect. 3.3.1], \cite{RainforthEtAl2024}; ~3) utilizing existing, easily computable estimates of conditional entropy or the MI. Since we aim to numerically maximize MI, which also depends on the design parameter $U$, the procedure for calculating MI will be called by the optimization solver millions of times and must therefore be sufficiently fast. Consequently, although MC methods provide good estimates of MI, they are of limited use here. The VLB methods require simultaneous optimization of the variational distribution with respect to its parameters and the signal $U$, \cite{HuanJagalurMarzouk2024}[Sect. 3.3.1 and 4.3.4], \cite{RainforthEtAl2024}. In addition, stochastic simulations are also used to compute the expected values in the VLB. As the goal of this article is to develop a simple design method that does not require hours of computation, we focus on the third option and use existing, easily computable lower-bounds of entropy or MI provided in \cite{Kolchinsky2017} or \cite{Huber2008}.

\section{Approximate solutions}
\label{sec:appr_solutions}
The optimization problem ~\eqref{eq:optimization_problem} becomes considerably more tractable when the parameter space $\vec{\Theta}$ is finite. Consequently, we first derive an approximate solution for a finite set of parameters and subsequently extend this result to obtain an approximate solution for the case where $\vec{\Theta}$ is an uncountable subset of $\mathbb{R}^{n_{\theta}}$.

\subsection{Finite parameter space}
\label{sec:section_3_1}
Let $\mat{\Theta}=\lbrace\vec{\theta}_1,...,\vec{\theta}_r\rbrace, \vec{\theta}_i \in \mathbb{R}^{n_{\theta}}$, $\vec{\theta}_i\neq\vec{\theta}_j$ and assume that $p_0$ is a discrete distribution of the form
\begin{equation}
p_0(\vec{\theta})=\sum\limits_{j=1}^{r}p_{0, j}\delta(\vec{\theta}-\vec{\theta}_j),
\label{eq:discrete_th_distr}
\end{equation}
where $p_{0,j}=P(\vec{\theta}=\vec{\theta}_j)$. Then, on the basis of \eqref{eq:pdfs_00} and \eqref{eq:discrete_th_distr}, the density of $\vec{Y}$ becomes a Gaussian mixture:
\begin{equation}
p(\vec{Y}|\vec{U})=\sum\limits_{j=1}^{r}p_{0, j}\mathcal{N}(\vec{Y},\vec{F}(\vec{\theta}_j, \vec{U}), \mat{S}(\vec{\theta}_j, \vec{U})).
\label{eq:discrete_evidence}
\end{equation}
The application of the Bayes rule gives the posterior:
\begin{equation}
p(\vec{\theta}_j|\vec{Y}, \vec{U})=\frac{p_{0, j}\mathcal{N}(\vec{Y},\vec{F}(\vec{\theta}_j, \vec{U}), \mat{S}(\vec{\theta}_j, \vec{U}))}{p(\vec{Y}|\vec{U})}.
\label{eq:discrete_posterior}
\end{equation}
The discrete counterpart of formulas \eqref{eq:entropies_00}-\eqref{eq:MI} takes the form:
\begin{align}
&H_{\vec{Y}}(\vec{U})=-\int p(\vec{Y}|\vec{U}) \ln p(\vec{Y}|\vec{U})d\vec{Y},\\
&H_{\vec{\theta}}=-\sum\limits_{j=1}^r p_{0, j}\ln p_{0, j},\\
&H_{\vec{Y}|\vec{\theta}}(\vec{U})=\tfrac{1}{2}\sum\limits_{j=1}^r p_{0, j}\ln \left((2\pi e)^{n_Y}|\mat{S}(\vec{\theta}_j, \vec{U})|\right),\\
&H_{\vec{\theta}|\vec{Y}}(\vec{U})=\int p(\vec{Y}|\vec{U})\left(-\sum\limits_{j=1}^r p(\vec{\theta}_j|\vec{Y}, \vec{U})\ln p({\vec{\theta}_j}|\vec{Y}, \vec{U})\right)d\vec{Y},
\label{eq:discrete_entropies}
\end{align}
\begin{equation}
I_{\vec{\theta};\vec{Y}}(\vec{U})=H_{\vec{\theta}}-H_{\vec{\theta}|\vec{Y}}(\vec{U})=H_{\vec{Y}}(\vec{U})-H_{\vec{Y}|\vec{\theta}}(\vec{U}).
\label{eq:MI_discrete}
\end{equation}
The direct computation of mutual information \eqref{eq:MI_discrete} remains difficult and, in many cases, intractable. \textbf{Hence, our central idea is to overcome this difficulty by replacing mutual information \eqref{eq:MI_discrete} with a computationally tractable and non-trivial lower bound.} In particular, we observe that \eqref{eq:discrete_evidence} is a finite Gaussian mixture. For such mixtures, one of the most effective lower bounds on $I_{\vec{\theta};\vec{Y}}$ is the inequality introduced in \cite{Kolchinsky2017}.

\begin{lem} (Information bounds, \cite{Kolchinsky2017}). For the Gaussian mixture \eqref{eq:discrete_evidence} with $p_{0,j}=P(\vec{\theta}=\vec{\theta}_j)$, the following inequalities holds:
\begin{equation}
I_l(\vec{U})\leqslant I_{\vec{\theta};\vec{Y}}(\vec{U}) \leqslant H_{\vec{\theta}},
\label{eq:inf_lb}
\end{equation}
where
\begin{align}
\label{eq:inf_lb_formulas_0}
& I_l(\vec{U})=-\sum\limits_{i=1}^rp_{0,i}\ln \left(\sum\limits_{j=1}^rp_{0,j} e^{-d_{i,j}(\vec{U})} \right),\\
\label{eq:inf_lb_formulas_dij}
& d_{i,j}(\vec{U})=\tfrac{1}{8}\vec{\Delta}_{i,j}^{\TT}\left(\tfrac{1}{2}(\mat{S}_i+\mat{S}_j)\right)^{-1}\vec{\Delta}_{i,j}+\tfrac{1}{2}\ln |\tfrac{1}{2}(\mat{S}_i+\mat{S}_j)|-\tfrac{1}{4}\ln\left(|\mat{S}_i||\mat{S}_j|\right),\\
& \vec{\Delta}_{i,j}=\vec{F}(\vec{\theta}_i, \vec{U})-\vec{F}(\vec{\theta}_j, \vec{U}), \mat{S}_i=\mat{S}(\vec{\theta}_i, \vec{U}), \mat{S}_j=\mat{S}(\vec{\theta}_j, \vec{U}).
\label{eq:inf_lb_formulas}
\end{align}
\label{lem:Lemma_1}
\end{lem}
Detailed proof is given in \cite{Kolchinsky_arxiv}[Sect. IIIb, inequality (11), and Sect. IV, formula (15) with $\alpha=0.5$] and also in \cite{Kolchinsky2017}.
Now, the approximate solution of \eqref{eq:optimization_problem} is given by  
\begin{equation}
\vec{U}^{\ast}=\mathrm{arg}\max_{\vec{U}\in\mathcal{U}_{ad}} I_l(\vec{U}).
\label{eq:optimization_problem_discrete}
\end{equation}
Since $\mathcal{U}_{ad}$ is compact and $I_l$ is smooth and bounded, the solution of \eqref{eq:optimization_problem_discrete} exists. We also note that in the case of two alternatives, that is, when $r=2$ in \eqref{eq:inf_lb_formulas_0}, we get  
\begin{equation}
    e^{-I_l(\vec{U})}=\left(p_{0,1}+p_{0,2}e^{-d_{1,2}(\vec{U})}\right)^{p_{0,1}}\left(p_{0,1}e^{-d_{1,2}(\vec{U})}+p_{0,2}\right)^{p_{0,2}}.
    \label{eq:eil_two_poss}
\end{equation}
Accordingly, the optimal signal in this case arises as the solution of a somewhat simplified optimization problem:
\begin{equation}
\max_{\vec{U}\in\mathcal{U}_{ad}} d_{1,2}(\vec{U}).
\label{eq:opt_problem_two_poss}
\end{equation}

If the function $\vec{F}$ in \eqref{eq:model_01} is affine with respect to $\vec{U}$ and the covariance $\mat{S}$ does not depend on $\vec{U}$, then it follows from Lemma \ref{lem:Lemma_1}, that $d_{1,2}$ is a positive (semi-) definite quadratic form with respect to $\vec{U}$. For constraints \eqref{eq:U_ad_1}, we thus obtain a convex quadratic programming problem. In the case of constraints \eqref{eq:U_ad}, one needs to find the minimum of $d_{1,2}$ on a closed ball in $\mathbb{R}^{n_U}$. This is also a convex problem, and it can be reduced to finding zeros of a scalar function \cite{NocedalWright2006}[Thm. 4.1, p. 70 and Sect. 4.3]. Furthermore, if $\vec{F}(\vec{\theta}_i,\vec{U})=\vec{F}_i\vec{U}, i=1,2$, and the constraints are defined by \eqref{eq:U_ad}, then the solution of \eqref{eq:opt_problem_two_poss} is the eigenvector of the matrix $\mat{Q}=(\vec{F}_1-\vec{F}_2)^T(\mat{S}_1+\mat{S}_2)^{-1}(\vec{F}_1-\vec{F}_2)$, corresponding to its largest eigenvalue (see \cite{Bania2019} [Sect. 2.1], for details).

\subsection{Infinite parameter space}
Let us assume that $\mat{\Theta}=\mathbb{R}^{n_{\theta}}$ and consider the Gaussian prior 
\begin{equation} 
p_0(\vec{\theta})=\cN(\vec{\theta},\vec{m}_{\theta},\mat{S}_{\theta}), \mat{S}_{\theta}>0.
\label{eq:gaussian_prior}
\end{equation}
Then the integral in \eqref{eq:pdfs_00} can be approximated with a finite Gaussian mixture
\begin{equation}
p(\vec{Y}|\vec{U})=\int p_0(\vec{\theta})\cN(\vec{Y},\vec{F}(\vec{\theta}, \vec{U}), \mat{S}(\vec{\theta}, \vec{U}))d\vec{\theta}\approx\sum\limits_{j=1}^{N_a}p_{0, j}\cN(\vec{Y},\vec{F}(\vec{\theta}_j, \vec{U}), \mat{S}(\vec{\theta}_j, \vec{U})),
\label{eq:finite_gauss_mix_appr}
\end{equation}
where $p_{0, j}\geqslant 0$ and $\sum\limits_{j=1}^{N_a} p_{0, j}=1$. 
The weights $p_{0, j}$ and the nodes $\vec{\theta}_j$ in \eqref{eq:finite_gauss_mix_appr} can be calculated by using the multidimensional Gauss-Hermite quadrature rule or any other suitable method. 

The Gauss-Hermite quadrature of order $p$ is exact for polynomials of degree at most $2p - 1$. The approximation error of the integral \eqref{eq:finite_gauss_mix_appr} using the Gauss-Hermite quadrature of order~$p$ 
depends on the $2p$-th derivatives of the functions $F$ and $S$. In the single-parameter case, with Gaussian prior $\cN(\theta, m_{\theta}, \sigma_{\theta})$, the error estimate is given by the formula:
\begin{equation}
    e\leqslant\sigma_{\theta}\frac{4^pC}{p!}\sup_{\theta,\vec{U},\vec{Y}}\frac{d^{2p}}{d\theta^{2p}}\cN(\vec{Y},\vec{F}(\theta, \vec{U}), \mat{S}(\theta, \vec{U})),
    \label{eq:gh_error_estimates}
\end{equation}
where $C$ is constant. The error tends to zero as $p\rightarrow\infty$ or $\sigma_{\theta}\rightarrow 0$. Therefore, the Gauss-Hermite approximation of the integral \eqref{eq:finite_gauss_mix_appr} is especially useful when the prior is narrow (small $\sigma_{\theta}$) or when the integrand, in the neighbourhood of the point $m_{\theta}$, can be well approximated by low-degree polynomials. To illustrate the method, we will show only a very simple second-order Gaussian quadrature rule with $2n_{\theta}$ points. 
 \begin{lem}
Approximate value of the integral $J(f)=\int \cN(\vec{\theta},\vec{m}_{\theta}, \mat{S}_{\theta})f(\vec{\theta})d\vec{\theta}$ is given by
 \begin{equation}
 J(f)\approx \tfrac{1}{2n_{\theta}}\sum\limits_{j=1}^{2n_{\theta}}f(\vec{\theta}_j),
 \label{eq:appr_gauss_integral}
 \end{equation}
 where
 \begin{equation}
\vec{\theta}_{2i-1}=\vec{m}_{\theta}-\mat{S}_{\theta}^{0.5}\sqrt{n_{\theta}}\vec{e}_i, \vec{\theta}_{2i}=\vec{m}_{\theta}+\mat{S}_{\theta}^{0.5}\sqrt{n_{\theta}}\vec{e}_i, i=1, ..., n_{\theta}
\label{eq:thetas}
\end{equation}
and $\vec{e}_i$ is $i^{th}$ basis vector in $\mathbb{R}^{n_{\theta}}$. If $f(\vec{\theta})=\tfrac{1}{2}\vec{\theta}^{\TT}\mat{A}\vec{\theta}+\vec{b}^{\TT}\vec{\theta}+c$, then the equality holds in \eqref{eq:appr_gauss_integral}.
\label{lem:Lemma_2}
\end{lem}
\begin{proof} Direct calculation. \end{proof}

An analogous method can be used for prior distributions defined on compact subsets of $\mathbb{R}^{n_{\theta}}$ (e.g. an $n$-dimensional hypercube), but the formulas for nodes and weights in \eqref{eq:finite_gauss_mix_appr} will then change. For example, if $\theta$ is a scalar parameter and the prior distribution is uniform, that is, $p_0(\theta) = U[a,b]$, then, using a second-order Gauss-Legendre quadrature, the approximate value of the integral $\int\limits_a^b p_0(\theta)f(\theta)d\theta$ is computed using the formula:
\begin{equation}
    \int\limits_a^b p_0(\theta)f(\theta)d\theta\approx p_{0,1}f(\theta_1)+p_{0,2}f(\theta_2),
    \label{eq:gauss_legendre_formula}
\end{equation}
where $p_{0,1}=p_{0,2}=0.5$ and
\begin{equation}
    \theta_1=\tfrac{1}{2}\left(a+b-\frac{b-a}{\sqrt{3}}\right), \theta_2=\tfrac{1}{2}\left(a+b+\frac{b-a}{\sqrt{3}}\right).
    \label{eq:gauss_legendre_nodes}
\end{equation}
The formula \eqref{eq:gauss_legendre_formula} is exact for polynomials of degree 3. The error estimate is analogous to \eqref{eq:gh_error_estimates} and tends to zero when $b-a\rightarrow 0$. More general multidimensional formulas, integration methods and error estimates are given in \cite{DavisRabinowitz1984}, \cite{Stroud1971}, \cite{Smolyak1963}.

The application of lemma \ref{lem:Lemma_2} to \eqref{eq:finite_gauss_mix_appr} gives $N_a=2n_{\theta}$, $p_{0, j}=(2n_{\theta})^{-1}$. Now, since \eqref{eq:finite_gauss_mix_appr} is approximated by a Gaussian mixture, the results of section \ref{sec:section_3_1} can be used. Based on \eqnsref{eq:inf_lb_formulas_0}{eq:inf_lb_formulas}, \eqref{eq:finite_gauss_mix_appr} and Lemma \ref{lem:Lemma_1}, the information lower bound takes the form: 
\begin{equation}
I_l(\vec{U})=-\frac{1}{2n_{\theta}}\sum\limits_{i=1}^{2n_{\theta}}\ln \left(\frac{1}{2n_{\theta}}\sum\limits_{j=1}^{2n_{\theta}}e^{-d_{i,j}(\vec{U})} \right),
\label{eq:inf_lb_simple}
\end{equation}
where $d_{i,j}$ and $\vec{\theta}_j$ are given by \eqnsref{eq:inf_lb_formulas_0}{eq:inf_lb_formulas} and \eqref{eq:thetas} or \eqref{eq:gauss_legendre_nodes}. The approximate solution of \eqref{eq:optimization_problem} can be found by maximizing \eqref{eq:inf_lb_simple} with constraints \eqref{eq:U_ad} or \eqref{eq:U_ad_1}.
%***********************************
\section{Bayesian input signal design in quasi-linear control systems}
\label{sec:quasi_linear_sys}
Consider the family of quasi-linear systems
\begin{equation}
\vec{x}_{k+1} =\mat{A}(\vec{\theta}, \vec{u}_k)\vec{x}_{k}+\mat{B}(\vec{\theta}, \vec{u}_k)+\mat{G}(\vec{\theta}, \vec{u}_k)\vec{w}_{k},
\label{eq:sys_dynamics}
\end{equation}
\begin{equation}
\vec{y}_{k} =\mat{C}\vec{x}_{k}+\vec{v}_k,
\label{eq:sys_observation}
\end{equation}
where $k=0,1, ..., N$, $N\geqslant 1$, $\vec{x}_{k}\in\mathbb{R}^{n}, \vec{y}_{k}\in\mathbb{R}^{n_y}, \vec{w}_{k}\in\mathbb{R}^{n_{w}}, \vec{v}_{k}\in\mathbb{R}^{n_y}$, $\vec{w}_{k}\sim\,\cN(0,\mathbb{I})$, $\vec{v}_{k}\sim\,\cN(0,\mat{S}_v), \mat{S}_v>0$. Variables $ \vec{x}_0, \vec{w}_{0}, ..., \vec{w}_{N-1}, \vec{v}_0, ..., \vec{v}_{N}$ are mutually independent. The initial state $\vec{x}_0$ is conditionally normal i.e. $p(\vec{x}_0|\vec{\theta})=\cN(\vec{x}_0,\vec{m}_0^-(\vec{\theta}),\mat{S}_0^-(\vec{\theta}))$, where $\vec{m}_0^-$, $\mat{S}_0^-$ are smooth and $\mat{S}_0^-(\vec{\theta})>0$, for all $\vec{\theta}\in\vec{\Theta}$. The joint prior distribution of the initial state $\vec{x}_0$ and the parameter $\vec{\theta}$ is given by $p_0(\vec{x}_0, \vec{\theta}) = p_0(\vec{\theta})\cN(\vec{x}_0, \vec{m}_0^-(\theta), \mat{S}_0^-(\vec{\theta}))$.
Let's define $\mat{A}_k=\mat{A}(\vec{\theta},\vec{u}_k)$, $\mat{B}_k=\mat{B}(\vec{\theta},\vec{u}_k)$, $\mat{G}_k=\mat{G}(\vec{\theta},\mat{u}_k)$. 
Then the solution of \eqref{eq:sys_dynamics} has the form:
\begin{align}
    \label{eq:sys_solution_0}
    &\vec{x}_0=\mathbb{I}\vec{x}_0,\\
    &\vec{x}_1=\mat{A}_0\vec{x}_0+\mat{B}_0+\mat{G}_0\vec{w}_0,\\
    &\vec{x}_2=\mat{A}_1\vec{x}_1+\mat{B}_1+\mat{G}_1\vec{w}_1=\mat{A}_1\mat{A}_0\vec{x}_0+\mat{A}_1\mat{B}_0+\mat{B}_1+\mat{A}_1\mat{G}_0\vec{w}_0+\mat{G}_1\vec{w}_1,\\
    &\vdots\\
    &\vec{x}_{N}=\vec{\Phi}(N,0)\vec{x}_0+\sum\limits_{j=0}^{N-1}\vec{\Phi}(N,j+1)\vec{B}_j+\sum\limits_{j=0}^{N-1}\vec{\Phi}(N,j+1)\mat{G}_j\vec{w}_j,
    \label{eq:sys_solution}
\end{align}
where $\vec{\Phi}(n,n)=\mathbb{I}$ and
\begin{equation}
    \vec{\Phi}(n,j)=\prod\limits_{i=1}^{n-j}\mat{A}_{n-i}, j<n.
    \label{eq:A_prod}
\end{equation}
Now, if we denote $\vec{X}=\mathrm{col}(\vec{x}_0,...,\vec{x}_N)$, $\vec{Y}=\mathrm{col}(\vec{y}_0,...,\vec{y}_N)$, $\vec{U}=\mathrm{col}(\vec{u}_0,...,\vec{u}_{N-1})$, $\vec{W}=\mathrm{col}(\vec{w}_0,...,\vec{w}_{N-1})$, $\vec{V}=\mathrm{col}(\vec{v}_0,...,\vec{v}_N)$, we can rewrite \eqnsref{eq:sys_solution_0}{eq:sys_solution} and \eqref{eq:sys_observation}, in matrix-vector form
\begin{equation}
\vec{X}=\mathcal{A}(\vec{\theta},\vec{U})\vec{x}_0+\mathcal{B}(\vec{\theta},\vec{U})+\mathcal{G}(\vec{\theta},\vec{U})\vec{W},
\label{eq:sys_solution_mtx}
\end{equation}
\begin{equation}
\vec{Y}=\mathcal{C}\vec{X}+\vec{V},
\label{eq:sys_observation_mtx}
\end{equation}
where the matrices $\mathcal{A}$, $\mathcal{B}$, $\mathcal{G}$, $\mathcal{C}=\mathbb{I}_{N+1}\otimes\mat{C}$ 
follow directly from \eqnsref{eq:sys_solution_0}{eq:sys_solution}, \eqref{eq:sys_observation}, and $\vec{W}\sim \cN(0,\mathbb{I}_{Nn_w})$, $\vec{V}\sim \cN(0,\mathbb{I}_{N+1}\otimes\mat{S}_v)$. Substituting \eqref{eq:sys_solution_mtx} into \eqref{eq:sys_observation_mtx} and taking into account that $p(\vec{x}_0|\vec{\theta})=\cN(\vec{x}_0,\vec{m}_0^-(\vec{\theta}),\mat{S}_0^-(\vec{\theta}))$, we get 
\begin{align}
\label{eq:model_lin_sys_0}
&\vec{Y}=\mathcal{C}\mathcal{A}(\vec{\theta},\vec{U})\vec{m}_0^-(\vec{\theta})+\mathcal{C}\mathcal{B}(\vec{\theta},\vec{U})+\vec{Z},\\
&\vec{Z}=\mathcal{C}\mathcal{A}(\vec{\theta},\vec{U})(\vec{x}_0-\vec{m}_0^-(\vec{\theta}))+\mathcal{C}\mathcal{G}(\vec{\theta},\vec{U})\vec{W}+\vec{V}.
\label{eq:model_lin_sys}
\end{align}
The conditional density of variable $\vec{Z}$ has the form $p(\vec{Z}|\vec{\theta})=\cN(\vec{Z}, 0, \mat{S}(\vec{\theta},\vec{U}))$, where the covariance matrix $\mat{S}$ is given by
\begin{equation}
\mat{S}(\vec{\theta},\vec{U})=\mathcal{C}(\mathcal{A}(\vec{\theta},\vec{U})\mat{S}_0^-(\vec{\theta})\mathcal{A}(\vec{\theta},\vec{U})^{\TT}+\mathcal{G}(\vec{\theta},\vec{U})\mathcal{G}(\vec{\theta},\vec{U})^{\TT})\mathcal{C}^{\TT}+\mathbb{I}_{N+1}\otimes\mat{S}_v.
\label{eq:total_cov_mtx}
\end{equation}
Finally, if we define 
\begin{equation}
    \vec{F}(\vec{\theta}, \vec{U})=\mathcal{C}\mathcal{A}(\vec{\theta},\vec{U})\vec{m}_0(\vec{\theta})+\mathcal{C}\mathcal{B}(\vec{\theta},\vec{U}),
    \label{eq:big_f_def_quasi_lin}
\end{equation}
we can rewrite \eqref{eq:model_lin_sys_0} in the form $\vec{Y}=\vec{F}(\vec{\theta}, \vec{U})+\vec{Z}$, which is exactly the model \eqref{eq:model_01}. To find the optimal input signal, we maximize one of the criteria \eqref{eq:inf_lb_formulas_0} or \eqref{eq:inf_lb_simple} with constraints \eqref{eq:U_ad} or \eqref{eq:U_ad_1}.
%***************************
%\subsection{Calculating the lower bound of mutual information}

With a large number of data (large $N$), calculating the inverse and determinant of a very large matrix
$\mat{S}(\vec{\theta}, \vec{U})$ in \eqref{eq:log_likelihood_0} and calculating the quantities $d_{i,j}(\vec{U})$ in \eqnsref{eq:inf_lb_formulas_0}{eq:inf_lb_formulas} is numerically ill-conditioned and requires special treatment. The algorithms below reduce the size of the matrices necessary to invert by a factor of $N+1$. 
\begin{lem}[Efficient computation of log-likelihood]
\label{lem:Lemma_3}
The following identities hold:
\begin{equation}
p(\vec{Y}|\vec{\theta}, \vec{U})=\prod\limits_{k=0}^N\cN(\vec{y}_k, \mat{C}\vec{m}_k^-(\vec{\theta}), \vec{\Sigma}_k(\vec{\theta})),\\
\label{eq:NYF_calc}
\end{equation}
\begin{equation}
    |\vec{Y}-\vec{F}(\vec{\theta}, \vec{U})|_{\mat{S}(\vec{\theta}, \vec{U})^{-1}}^2=\sum\limits_{k=0}^N|\vec{y}_k-\mat{C}\vec{m}_k^-(\vec{\theta})|_{\vec{\Sigma}_k^{-1}(\vec{\theta})}^2,
    \label{eq:norm_Y_F_calc}
\end{equation}
\begin{equation}
    |\mat{S}(\vec{\theta}, \vec{U})|=\prod\limits_{k=0}^N|\vec{\Sigma}_k(\vec{\theta})|,
    \label{eq:det_S_calc}
\end{equation}
\begin{equation}
\mathcal{L}(\vec{\theta},\vec{Y},\vec{U})=\tfrac{1}{2}\sum\limits_{k=0}^{N}\left(|\vec{y_k}-\mat{C}\vec{m}_k^-(\vec{\theta})|_{\vec{\Sigma}_k^{-1}(\vec{\theta})}^2+\ln|\vec{\Sigma}_k(\theta)|\right)-\ln p_0(\vec{\theta}),
\label{eq:log_likelihood_recursive}
\end{equation}
where $\mathcal{L}$ is given by \eqref{eq:log_likelihood_0} and $\vec{m}_k^-$, $\vec{\Sigma}_k$, are calculated recursively by the Kalman filter
\begin{align}
\label{eq:pem_kalman_0}
&\vec{\Sigma}_k(\vec{\theta})=\mat{S}_v+\mat{C}\mat{S}_k^-(\vec{\theta})\mat{C}^{\TT},\\
&\mat{L}_k(\vec{\theta})=\mat{S}_k^-(\vec{\theta})\mat{C}^{\TT}\vec{\Sigma}_k^{-1}(\vec{\theta}),\\
&\vec{m}_k(\vec{\theta})=\vec{m}_k^-(\vec{\theta})+\mat{L}_k(\vec{\theta})(\vec{y}_k-\mat{C}\vec{m}_k^-(\vec{\theta})),\\
&\mat{S}_k(\vec{\theta})=\mat{S}_k^-(\vec{\theta})-\mat{L}_k(\vec{\theta})\vec{\Sigma}_k(\vec{\theta})\mat{L}_k(\vec{\theta})^{\TT},\\
&\vec{m}_{k+1}^-(\vec{\theta})=\mat{A}_{k}\vec{m}_{k}(\vec{\theta})+\mat{B}_{k},\\
&\mat{S}_{k+1}^-(\vec{\theta})=\mat{A}_{k}\mat{S}_{k}(\vec{\theta})\mat{A}_{k}^{\TT}+\mat{G}_{k}\mat{G}_{k}^{\TT}, k=0, 1 ..., N,
\label{eq:pem_kalman}
\end{align}
with initial conditions $\vec{m}_0^-(\vec{\theta})$, $\mat{S}_{0}^-(\vec{\theta})$.
\end{lem}
The proof is given in Appendix \ref{app:appendix_A}. The \eqnsref{eq:pem_kalman_0}{eq:pem_kalman}, are, in fact, a family of discrete-time Kalman filters indexed by $\vec{\theta}$. The first four formulas describe the correction step. The prediction step is given by the last two equations. The matrix $\mat{L}_k$ is the Kalman gain and $\vec{\Sigma}_k$ is the covariance matrix of the output prediction error $\vec{\epsilon}_k=\vec{y}_k-\mat{C}\vec{m}_k^-$. 
\begin{lem}[Efficient computation of $d_{i,j}$ in \eqref{eq:inf_lb_formulas_dij}]
\label{lem:lemma_4}
Let us define
\begin{equation}
\tilde{\mat{A}}_{k}=\begin{bmatrix}
    \mat{A}(\vec{\theta}_i,\vec{u}_k) & 0\\
    0&\mat{A}(\vec{\theta}_j,\vec{u}_k)
\end{bmatrix},
\tilde{\mat{B}}_{k}=\begin{bmatrix}
    \mat{B}(\vec{\theta}_i,\vec{u}_k)\\
    \mat{B}(\vec{\theta}_j,\vec{u}_k)
\end{bmatrix},
\label{eq:system_2n_ab}
\end{equation}
\begin{equation}
\tilde{\mat{G}}_{k}=\begin{bmatrix}
    \mat{G}(\vec{\theta}_i,\vec{u}_k) & 0\\
    0&\mat{G}(\vec{\theta}_j,\vec{u}_k)
\end{bmatrix},
\tilde{\mat{C}}=\frac{1}{\sqrt{2}}\begin{bmatrix}
    \mat{C} & -\mat{C}
\end{bmatrix}
\label{eq:system_2n_gc}
\end{equation}
and let
\begin{align}
\label{eq:kalman_ij_0}
&\tilde{\vec{\Sigma}}_k=\mat{S}_v+\tilde{\mat{C}}\tilde{\mat{S}}_k^-\tilde{\mat{C}}^{\TT},\\
&\tilde{\mat{L}}_k=\tilde{\mat{S}}_k^-\tilde{\mat{C}}^{\TT}\tilde{\vec{\Sigma}}_k^{-1},\\
&\tilde{\mat{S}}_k=\tilde{\mat{S}}_k^--\tilde{\mat{L}}_k\tilde{\vec{\Sigma}}_k\tilde{\mat{L}}_k^{\TT},\\
&\tilde{\vec{m}}_{k+1}^-=\tilde{\mat{A}}_{k}(\mathbb{I}-\tilde{\mat{L}}_k\tilde{\mat{C}})\tilde{\vec{m}}_{k}^-+\tilde{\mat{B}}_{k},\\
&\tilde{\mat{S}}_{k+1}^-=\tilde{\mat{A}}_{k}\tilde{\mat{S}}_{k}\tilde{\mat{A}}_{k}^{\TT}+\tilde{\mat{G}}_{k}\tilde{\mat{G}}_{k}^{\TT}, k=0, 1, ..., N,
\label{eq:kalman_ij_1}
\end{align}
with initial conditions
\begin{equation}
    \tilde{\vec{m}}_0^-=\begin{bmatrix}
        \vec{m}_0^-(\vec{\theta}_i)\\
        \vec{m}_0^-(\vec{\theta}_j)
    \end{bmatrix}, 
    \tilde{\mat{S}}_0^-=\begin{bmatrix}
        \mat{S}_0^-(\vec{\theta}_i)&0\\
        0&\mat{S}_0^-(\vec{\theta}_j)
    \end{bmatrix}.
    \label{eq:init_cond_d_ij}
\end{equation}
Then the quantity $d_{i,j}(\vec{U})$ in formula \eqref{eq:inf_lb_formulas_0} is given by:
\begin{equation}
    d_{i,j}(\vec{U})=\tfrac{1}{4}\sum\limits_{k=0}^N|\tilde{\vec{C}}\tilde{\vec{m}}_k^-|_{\tilde{\vec{\Sigma}}_k^{-1}}^2+\tfrac{1}{2}\sum\limits_{k=0}^N\ln|\tilde{\vec{\Sigma}}_k|-\tfrac{1}{4}\ln\left(|\mat{S}_i||\mat{S}_j|\right),
    \label{eq:d_ij_calc}
\end{equation}
where $|\mat{S}_i|=|\mat{S}(\vec{\theta_i,\vec{U}})|$, $|\mat{S}_j|=|\mat{S}(\vec{\theta_j,\vec{U}})|$, are calculated accordingly to Lemma \ref{lem:Lemma_3}, \eqnref{eq:det_S_calc}. 
\end{lem}
The proof is given in Appendix \ref{app:appendix_A}. Let us observe that, instead of calculating the inverse and determinant of the large matrices $\mat{S}_i$, $\mat{S}_j$, $\tfrac{1}{2}(\mat{S}_i+\mat{S}_j)$, of dimension $(N+1)n_y$, we only need to calculate the determinants and inverses of the much smaller matrices $\vec{\Sigma}_k$, $\tilde{\vec{\Sigma}}_k$, whose dimension is $n_y$, which is usually a small number.  

%******************
\section{Comparison with classical methods of input signal design}
\label{sec:classical_methods}
Classical methods for input signal design in system identification are primarily concerned with LTI state space or transfer function models (such as ARMAX) and are usually based on maximizing some functions of error covariance or Fisher information matrix. For the prediction error method (PEM) estimator, the asymptotic form of the error covariance matrix (or Fisher information) is well known, both in the time and frequency domains. In the time domain, the solution corresponds to a specific input signal, whereas in the frequency domain, the solution yields the optimal power spectral density of the input signal. Below, we provide a brief overview of these methods, following the methodology presented in \cite{GoodwinPayne1977}, \cite{Ljung1999}[Ch.9. Sect. 9.3, 9.4]  and \cite{Pronzato2008}[Sect. 6.1].

Consider the LTI, SISO system
\begin{align}
\label{eq:sys_classic}
    &\vec{x}_{k+1} = \mat{A}(\vec{\theta})\vec{x}_k + \mat{B}(\vec{\theta})u_k + \mat{G}(\vec{\theta})\vec{w}_k, \\
    &y_k = \vec{C}\vec{x}_k + v_k,
    \label{eq:sys_classic_1}
\end{align}
under the assumptions stated in Section~\ref{sec:quasi_linear_sys}. System \eqref{eq:sys_classic}, \eqref{eq:sys_classic_1} is equivalent to the transfer function model
\begin{equation}
    y_k=G(\vec{\theta},z)u_k+H(\vec{\theta},z)e_k,
    \label{eq:dt_tranfer_function_model}
\end{equation}
where $e_k\sim\cN(0,\sigma_e^2)$ is a sequence of mutually independent Gaussian variables. The filters $G$ and $H$ are determined by the formulas
\begin{equation}
    G(\vec{\theta},z)=\mat{C}(z\mathbb{I}-\mat{A}(\vec{\theta}))^{-1}\mat{B}(\vec{\theta}), H(\vec{\theta},z)=1+\mat{C}(z\mathbb{I}-\mat{A}(\vec{\theta}))^{-1}\mat{K}(\vec{\theta}),
    \label{eq:GH_filters_def}
\end{equation}
where the Kalman gain $\mat{K}(\vec{\theta})$ is given by  
\begin{equation}
\mat{K}(\vec{\theta})=\mat{A}(\vec{\theta})\mat{S}(\vec{\theta})\mat{C}^{\TT}(\mat{C}\mat{S}(\vec{\theta})\mat{C}^{\TT}+\sigma_v^2)^{-1},
\label{eq:kalman_gain_gh_filters}
\end{equation}
with a non-negative matrix $\mat{S}$ being a solution of the Riccati equation (cf. \cite{Ljung1999})
\begin{equation}
\mat{S}=\mat{A}\mat{S}\mat{A}^{\TT}+\mat{G}\mat{G}^{\TT}-\mat{A}\mat{S}\mat{C}^{\TT}(\mat{C}\mat{S}\mat{C}^{\TT}+\sigma_v^2)^{-1}\mat{C}\mat{S}\mat{A}^{\TT}.
\label{eq:riccati_eq_for_h}
\end{equation}
The prediction errors are given by the recurrence  
\begin{equation}
    \epsilon_k(\vec{\theta},\vec{Y},\vec{U})=H^{-1}(\vec{\theta},z)\left(y_k-G(\vec{\theta},z)u_k\right).
    \label{eq:prediction errors}
\end{equation}
The cost function used in the Prediction Error Method (PEM) is expressed as:
\begin{equation}
    V(\vec{\theta},\vec{Y}, \vec{U})=\frac{1}{2N\sigma_e^2}\sum\limits_{k=1}^N\epsilon_k^2(\vec{\theta},\vec{U}).
    \label{eq:pem_criterion}
\end{equation}
Minimization of \eqref{eq:pem_criterion} with reference to $\vec{\theta}$ yields the PEM estimator
\begin{equation}
    \hat{\vec{\theta}}(\vec{Y},\vec{U})=\mathrm{arg}\min_{\vec{\theta}\in\vec{\Theta}} V(\vec{\theta},\vec{Y},\vec{U}).
    \label{eq:pem_estimator}
\end{equation}
The above estimator, under rather weak identifiability conditions, is consistent, asymptotically normal, and efficient, i.e. it achieves the Cram\'{e}r-Rao lower bound. Following the reasoning presented in \cite{GoodwinPayne1977}, \cite{Ljung1999}[Ch. 9, Sect. 9.3 and 9.4] or \cite{Pronzato2008}[Sect. 6.1], we divide the parameter vector into two groups related to the parameters appearing in $G$ and $H$, that is, $\vec{\theta}=\mrm{col}(\vec{\theta_H}, \vec{\theta_G})$. The sensitivity of $\epsilon_k$ to changes in $\vec{\theta}_G$ is calculated recursively according to the following equations:
\begin{equation}
    \vec{\psi}_k(\vec{\theta},\vec{U})=H^{-1}(\vec{\theta},z)\left(\nabla_{\vec{\theta}_G}G(\vec{\theta},z)\right)u_k=F_z(\vec{\theta}, z)u_k,
    \label{eq:pem_sensitivity}
\end{equation}
where $\nabla_{\vec{\theta}_G}$ means differentiation only with respect to the parameters that occur in $G$. The information matrix, which is also the inverse of the error covariance $\mat{P}_{\vec{\theta}_G}$, is given by:  
\begin{equation}
    \mat{M}(\vec{\theta}, \vec{U})=\mat{P}_{\vec{\theta}_G}^{-1}(\vec{\theta}, \vec{U})=\mat{R}_e(\vec{\theta})+\frac{1}{N\sigma_e^2}\sum\limits_{k=1}^N\vec{\psi}_k(\vec{\theta},\vec{U})\vec{\psi}_k(\vec{\theta},\vec{U})^{\TT},
    \label{eq:pem_fisher_mtx}
\end{equation}
where $\mat{R}_e$ does not depend on $\vec{U}$. Using the D-optimal criterion, the optimal signal is given by maximization of $\det\left(\mat{M}(\vec{\theta}^0,\vec{U})\right)$, where $\vec{\theta}^0$ is the true value of the parameter. Since $\vec{\theta}^0$ is unknown, one can use the prior distribution and maximize the average D-optimal criterion:
\begin{equation}
    Q(\vec{U})=\E_{p_0(\vec{\theta})}\det\left(\mat{M}(\vec{\theta},\vec{U})\right),
    \label{eq:av_d_optimal_criterion}
\end{equation}
with constraints \eqref{eq:U_ad} or \eqref{eq:U_ad_1}. The asymptotic error covariance can also be expressed in terms of the power spectral density of the input signal $u_k$. Let $\Phi_u$ denote the spectral density of $u_k$. As it was shown in \cite{Ljung1999}[p. 291] or \cite{Pronzato2008} [Sect. 6.2], we have:
\begin{equation}
    \mat{M}(\Phi_u, \vec{\theta}^0)=\mat{P}_{\vec{\theta}_G}^{-1}(\Phi_u, \vec{\theta}^0)=\frac{N}{2\pi\sigma_e^2}\int\limits_{-\pi}^{\pi}F_z(e^{i\omega},\vec{\theta}^0)F_z(e^{-i\omega}, \vec{\theta}^0)^{\TT}\Phi_u(\omega)d\omega+\mat{R}_e(\vec{\theta}^0),
    \label{eq:pem_fisher_spectral}
\end{equation}
where $F_z$ is defined by \eqref{eq:pem_sensitivity} and the term $\mat{R}_e$ in \eqref{eq:pem_fisher_spectral} does not depend on $\Phi_u$. Similarly as before, the parameter-averaged determinant of the matrix $\mat{M}$ is maximized with respect to $\Phi_u$, subject to the signal power and frequency constraints. Typically, the spectrum $\Phi_u$ is parametrized by a finite number of coefficients $c_k$, so that the resulting optimization problem is convex; see \cite{Jansson2004} for details. After performing spectral factorization of $\Phi_u$, a filter is obtained, whose input is white noise and whose output yields the optimal signal $u_k$. This has been implemented in the MOOSE-2 solver \cite{MOOSE2}. Unfortunately, MOOSE-2 does not allow for averaging over prior and involves \textit{true} unknown value of the parameter.

Numerous variants of the aforementioned methods can be found in the literature. For example, instead of the D-optimality criterion, one may also consider maximizing $\mathrm{tr}(\mat{M})$ or $\lambda_{\min}(\mat{M})$. However, the vast majority of methods are based on the principles stated above (see, e.g., \cite{Pronzato2008}), that is, maximization of some functions of the Fisher information matrix. Finally, we note that the above methods employ the classical optimality criterion, averaged only over the prior distribution. Consequently, they are not fully Bayesian and, following the terminology of \cite{Ryan2015}, should rather be referred to as \textit{pseudo-Bayesian} methods.

\section{Examples of input signal design}
In the following, we present four examples of optimal input signal design using both Bayesian and classical methods. 
Examples 1-3 are classical in nature and concern time-invariant linear stochastic systems. Examples 1 and 2 are elementary, while Example~3, taken from~\cite{JimenezMatinez2018}, addresses the design of a control signal for a paradigmatic model of the atomic sensor. The sensor is modelled as a harmonic oscillator with the natural frequency being the parameter of interest. In examples 1-3, the Bayesian approach is compared with classical methods. Maximization of the spectral criterion \eqref{eq:pem_fisher_spectral} was performed using the MOOSE-2 solver \cite{MOOSE2}, evaluated at $\vec{\theta}=\vec{m}_{\theta}$ with default parameters, that is, the input spectrum was FIR-type with 20 lags and the spectrum power constraint was set to 1 (prob.spectrum.signal.power.ub=1). There were no additional constraints on the shape of the spectrum.

Example~4, adapted from~\cite{Truullinou2021}, is more advanced and considers the design of the pump laser control signal in an optically pumped magnetometer. The magnetometer is modelled as a quasi-linear stochastic system, where the matrices $\mat{A}$, $\mat{B}$, and $\mat{G}$ depend nonlinearly on the control signal $u$. For this system, classical methods cannot be applied. Therefore, estimation errors are compared with the Information-Theoretic Lower Bound (ITB) provided in Theorem \ref{thm:theorem_1} and with the errors obtained by using an appropriately selected harmonic input signal.

In all examples, the errors were computed using the Monte Carlo method. The parameter $\vec{\theta}$ and the initial conditions $\vec{x}_0$ were sampled from the prior distribution $p_0(\vec{x}_0, \vec{\theta})$. Observations $y_0, \dots, y_N$ corresponding to the sampled parameters and initial conditions were then generated, and the error of the MAP estimator was calculated. This error was subsequently averaged over many repetitions of the procedure.

\subsection{Elementary example}
\label{sec:example_01}
We begin with a very simple first-order system  
\begin{align}
    &x_{k+1}=\theta_1x_k+\theta_2u_k+gw_k,\\
    &y_k=x_k+\sigma_vv_k, k=0, 1, ..., N,
    \label{eq:sys_example_01}
\end{align}
with $\sigma_v=0.1$, $g=0.01$. The parameter vector $\vec{\theta}=[\theta_1,\theta_2]^{\TT}$, has a prior distribution $p_0(\vec{\theta})=\cN(\vec{\theta},\vec{m}_{\theta}, \mat{S}_{\theta})$, where $\vec{m}_{\theta}=\begin{bmatrix} 0.8&0.2 \end{bmatrix}^{\TT}$, $\mat{S}_{\theta}=10^{-3}\mathbb{I}$. As assumed in section \ref{sec:quasi_linear_sys}, the initial condition $x_0$ is conditionally Gaussian, that is, $p(x_0|\vec{\theta})=\cN(x_0,m_0^-(\vec{\theta}),s_0^-(\vec{\theta}))$, with $m_0^-(\vec{\theta})=0$, $s_0^-(\vec{\theta})=0.01$. The length of the signal $N=100$ and the set of admissible signals is given by \eqref{eq:U_ad} with $\tilde{\vec{U}}=0$, that is, the norm of the signal cannot be greater than $\varrho$. 
To minimize the averaged D-optimal criterion \eqref{eq:av_d_optimal_criterion}, we need to calculate the sensitivity of the prediction error. The sensitivity equations \eqref{eq:pem_sensitivity} now take the form
\begin{align}
    &(1+(K(\theta_1)-\theta_1)z^{-1})(1-\theta_1z^{-1})\psi_{1,k}=\theta_2z^{-2}u_k,\\
    &(1+(K(\theta_1)-\theta_1)z^{-1})\psi_{2,k}=z^{-1}u_k,
    \label{eq:sensitivity_example_01}
\end{align}
where the Kalman gain $K(\theta_1)$ is given by \eqref{eq:kalman_gain_gh_filters}, \eqref{eq:riccati_eq_for_h} with $\mat{A}=\theta_1$, $\mat{G}=g$, $\mat{C}=1$.

\begin{figure}[H]
\centering
%\isPreprints{\centering}{} % Only used for preprints
\includegraphics[width=13 cm]{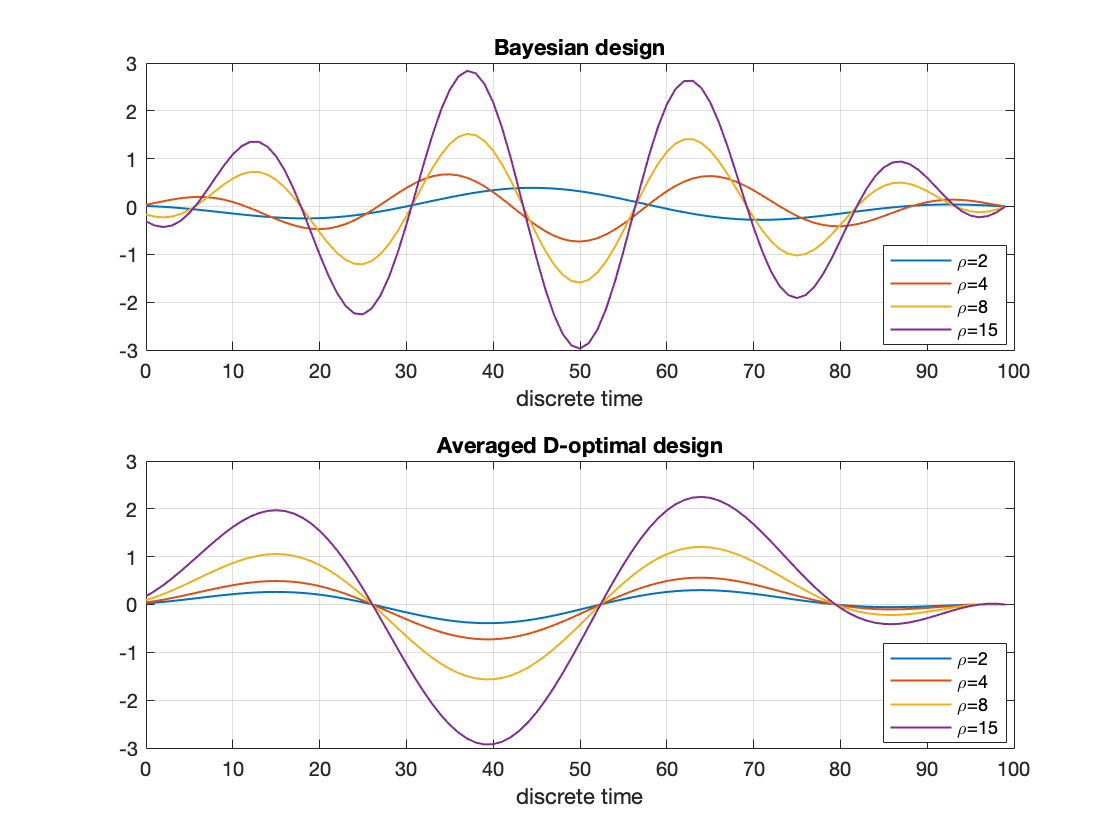}
\caption{Optimal input signals resulting from the maximization of the Bayesian criterion \eqref{eq:inf_lb_simple} (top) and of the averaged D-optimal criterion \eqref{eq:av_d_optimal_criterion} (bottom), shown for several values of $\varrho$.
\label{fig:example_01_f1}}
\end{figure}

\begin{figure}[H]
\centering
%\isPreprints{\centering}{} % Only used for preprints
\includegraphics[width=13 cm]{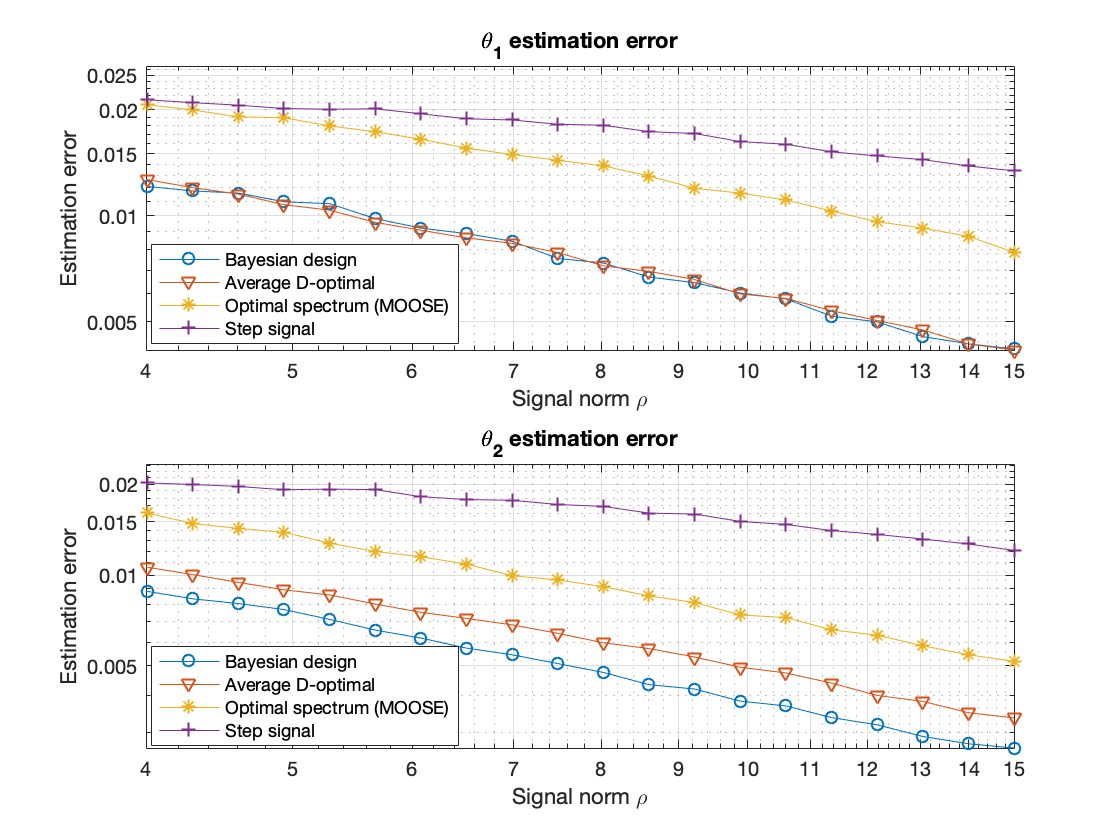}
\caption{Mean estimation errors of the parameters $\theta_1$ and $\theta_2$, obtained using the MAP estimator \eqref{eq:map_estimator}, as functions of the maximum admissible signal norm $\varrho$. The results are based on a Monte Carlo simulation with 3000 repetitions. The constant (step) signal and the MOOSE signal were always assigned a norm equal to $\varrho$.
	\label{fig:example_01_f2}}
\end{figure}

The optimal input signals were designed by maximizing the Bayesian criterion \eqref{eq:inf_lb_simple}, the averaged D-optimal criterion \eqref{eq:av_d_optimal_criterion}, and the spectral criterion \eqref{eq:pem_fisher_spectral}, subject to the constraint \eqref{eq:U_ad} with $\tilde{\vec{U}}=0$. The optimal signals and the corresponding estimation errors of $\theta_1$ and $\theta_2$ are shown in Figures \ref{fig:example_01_f1} and \ref{fig:example_01_f2}. In Fig. \ref{fig:example_01_f2} we also calculate the estimation errors for the constant (step) signal, which is certainly not optimal. The constant signal and the MOOSE signal were always assigned a norm equal to $\varrho$.   

\subsection{Example with a non-Gaussian prior distribution}
Consider the following system:
\begin{equation}
    d\vec{x}=\left(\mat{A}_c(\theta)+\mat{B}_c(\theta)u\right)dt+\mat{G}_c(\theta)d\vec{w},
    \label{eq:ct_sys_ex_02a}
\end{equation}
\begin{equation}
    y_k=\mat{C}\vec{x}_k+s_vv_k,
\end{equation}
where
\begin{equation}
    \mat{A}_c(\theta)=\begin{bmatrix}
        0&1\\
        0&-\theta
    \end{bmatrix}, 
    \mat{B}_c=\begin{bmatrix}
        0\\\theta
    \end{bmatrix},
    \mat{G}_c(\theta)=\begin{bmatrix}
        0\\\sqrt{d_c\theta}
    \end{bmatrix}, 
    \mat{C}=\begin{bmatrix}
        1&0
    \end{bmatrix},
    \label{eq:sys_matrices_ex_02_ct}
\end{equation}
$\vec{x}_k=\vec{x}(t_k)$, $t_k=k\Delta$, $\Delta=0.05\cdot 10^{-3}$, $d_c=0.01$, $s_v=0.1$. This system can be considered as controlled Brownian motion or DC motor with stochastic disturbances. The parameter $\theta$ is the unknown damping rate of the system. Assuming that $u(t)=u_k, t\in [t_k, t_{k+1}]$, the discrete-time system corresponding to \eqref{eq:ct_sys_ex_02a} has the form
\begin{equation}
    \vec{x}_{k+1}=\mat{A}(\theta)\vec{x}_k+\mat{B}(\theta)u_k+\mat{G}(\theta)w_k,
\end{equation}
where, according to the procedure given in Appendix \ref{app:appendix_C}:
\begin{equation}
\label{eq:sys_matrices_ex_02a}
\begin{split}
    &\mat{A}(\theta)=\begin{bmatrix}
        1 & \frac{1-e^{-\theta\Delta}}{\theta}\\
        0 & e^{-\theta\Delta}
    \end{bmatrix}, 
    \mat{B}(\theta)=\begin{bmatrix}
        \Delta-\frac{1-e^{-\theta\Delta}}{\theta}\\
        1-e^{-\theta\Delta}
    \end{bmatrix},\\
        &\mat{G}(\theta)=\sqrt{d_c\theta}\begin{bmatrix}
        \sqrt{D_{1,1}(\theta)} & 0\\
        \frac{D_{1,2}(\theta)}{\sqrt{D_{1,1}(\theta)}} & \sqrt{\frac{D_{1,1}(\theta)D_{2,2}(\theta)-D_{1,2}(\theta)^2}{D_{1,1}(\theta)}}
    \end{bmatrix},
    \end{split}
\end{equation}
where
\begin{align}
    & D_{1,1}(\theta)=\frac{4e^{-\theta\Delta}-e^{-2\theta\Delta}+2\theta\Delta-3}{2\theta^3},\\
    & D_{1,2}(\theta)=\frac{1-2e^{-\theta\Delta}+e^{-2\theta\Delta}}{2\theta^2}, D_{2,2}(\theta)=\frac{1-e^{-2\theta\Delta}}{2\theta}.
    \label{eq:d_matrix_ex_02a}
\end{align}
The initial condition is Gaussian with $\vec{m}_0^-=0$, $\mat{S}_0^-=\mrm{diag}[0.001, 0.005]$. Unlike in the previous example, here \textbf{we assume that the prior distribution of $\theta$ is uniform}, that is, $p_0(\theta)=U[a,b]$ with $a=0.05$, $b=2$. Following the Gauss-Legendre formula \eqref{eq:gauss_legendre_nodes}, we get $\theta_1=0.5(a+b-(b-a)/\sqrt{3})\approx 0.462$, $\theta_2=0.5(a+b+(b-a)/\sqrt{3})\approx 1.588$, $p_{0,1}=p_{0,2}=0.5$. Thus, $r=2$ in \eqref{eq:inf_lb_formulas_0} and, accordingly to \eqref{eq:eil_two_poss}, the Bayesian optimal signal is a solution of the simplified and convex optimization problem \eqref{eq:opt_problem_two_poss} with $d_{1,2}$ defined by Lemmas \ref{lem:Lemma_3} and \ref{lem:lemma_4}. Moreover, since the matrices in \eqref{eq:sys_matrices_ex_02a}, do not depend on $u_k$, the last two terms in \eqref{eq:d_ij_calc} can be omitted. The set of admissible signals is given by \eqref{eq:U_ad} with $\tilde{\vec{U}}=0$, that is, the signal norm cannot be greater than $\varrho$.

In order to employ the classical methods described in Sect. \ref{sec:classical_methods}, it is necessary to first evaluate the sensitivity of the prediction error. The transfer functions $G$ and $H$ in \eqref{eq:dt_tranfer_function_model} have the form
\begin{equation}
    G(\theta, z)=\frac{B(\theta,z)}{A(\theta,z)}z^{-1}, H(\theta, z)=\frac{C(\theta,z)}{A(\theta,z)},
    \label{eq:ex_02_gh_1}
\end{equation}
where
\begin{align}
 &A(\theta,z)=1-(1+e^{\theta\Delta}))z^{-1}+e^{-\theta\Delta}z^{-2},\\ 
 &B(\theta,z)=\Delta-\frac{1-e^{-\theta\Delta}}{\theta}+\left(\frac{1}{\theta}-\left(\frac{1}{\theta}+\Delta\right)e^{-\theta\Delta}\right)z^{-1},\\
 &C(\theta,z)=1+\left(\mat{K}_1(\theta)-\left(1-e^{-\theta\Delta}\right)\right)z^{-1}+\left(\frac{1-e^{-\theta\Delta}}{\theta}\mat{K}_2(\theta)+(1-\mat{K}_1(\theta))e^{-\theta\Delta}\right)z^{-2},
\label{eq:ex_02_polynomials} 
\end{align}
and the Kalman gain $\mat{K}$ is given by \eqref{eq:kalman_gain_gh_filters}. Since we only have one parameter, the sensitivity $\psi_k$ is a number, and the sensitivity equation \eqref{eq:pem_sensitivity} now takes the form
\begin{equation}
A(\theta,z)C(\theta,z)\psi_k(\theta, \vec{U})=\left(\frac{\partial B(\theta,z)}{\partial\theta}A(\theta,z)-B(\theta,z)\frac{\partial A(\theta,z)}{\partial\theta}\right)z^{-1}u_k.
\label{eq:ex_02_sensitivity}
\end{equation}
The D-optimal signal is then obtained by maximization of the averaged D-optimal criterion
\begin{equation}
    Q(\vec{U})=\E_{p_0(\theta)}\left[ \frac{1}{N}\sum\limits_{k=1}^N\psi_k^2(\theta, \vec{U})\right],
    \label{eq:av_d_optimal_ex_02}
\end{equation}
with constraints \eqref{eq:U_ad}.

\begin{figure}[H]
\centering
%\isPreprints{\centering}{} % Only used for preprints
\includegraphics[width=13 cm]{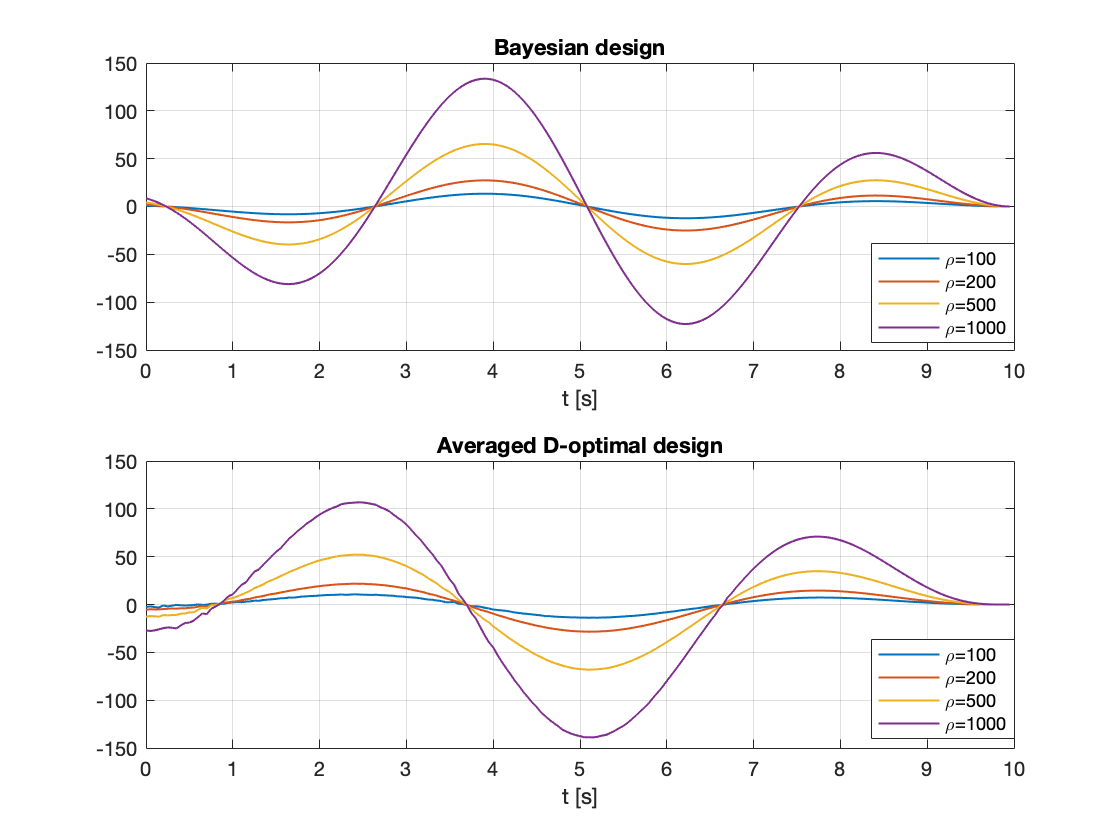}
\caption{Optimal input signals (left) and corresponding system outputs (right) obtained by maximizing the Bayesian criterion \eqref{eq:opt_problem_two_poss} (top) and the averaged D-optimal criterion \eqref{eq:av_d_optimal_ex_02} (bottom) subject to the constraint \eqref{eq:U_ad} with $\tilde{\vec{U}}=0$.
\label{fig:example_02a_f1}}
\end{figure}

\begin{figure}[H]
\centering
%\isPreprints{\centering}{} % Only used for preprints
\includegraphics[width=13 cm]{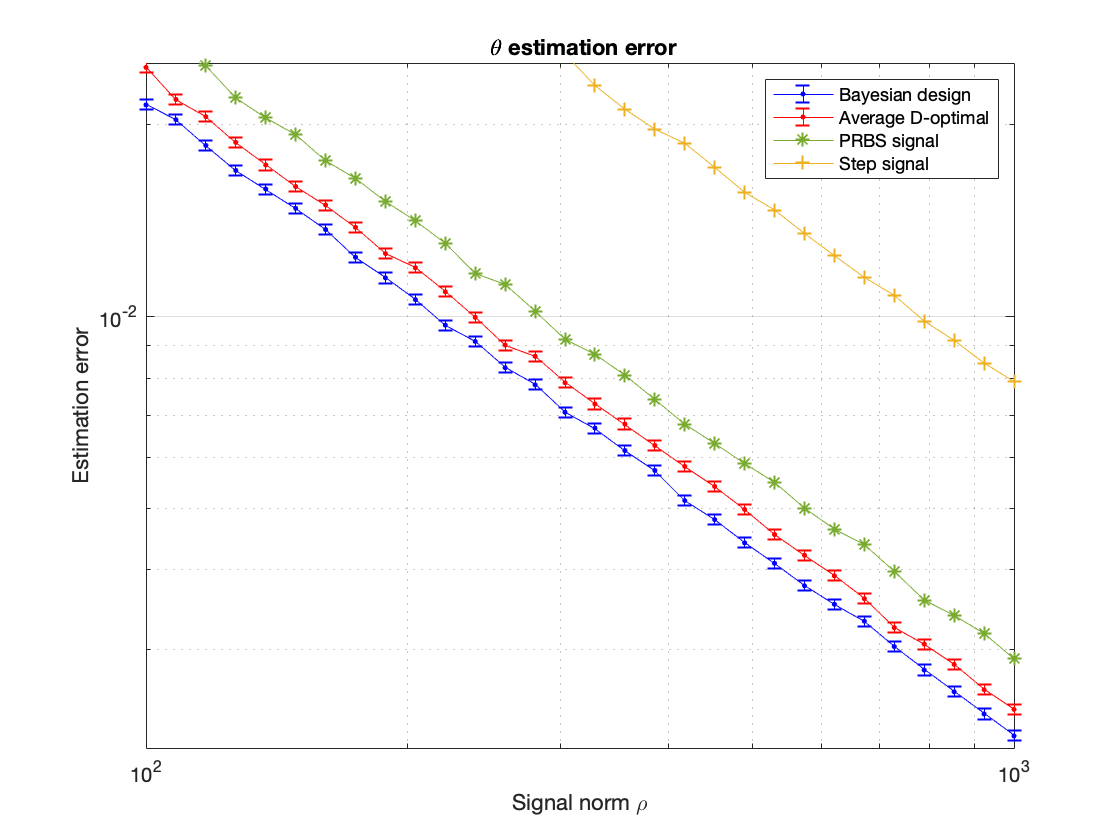}
\caption{Mean estimation errors of the parameter $\theta$, obtained using the MAP estimator \eqref{eq:map_estimator}, as functions of the maximum admissible signal norm $\varrho$, for different signals. The results are based on a Monte Carlo simulation with 6000 repetitions. Error bars show that the difference between the D-optimal and Bayesian methods is statistically significant.
\label{fig:example_02a_f2}}
\end{figure}

The optimal input signals were designed by maximizing the Bayesian criterion \eqref{eq:inf_lb_simple} and the averaged D-optimal criterion \eqref{eq:av_d_optimal_ex_02}, subject to the constraint \eqref{eq:U_ad} with $\tilde{\vec{U}}=0$. The results are presented in Figures \ref{fig:example_02a_f1} and \ref{fig:example_02a_f2}. Figure \ref{fig:example_02a_f2} also shows the estimation error for the step (i.e. constant) signal and the PRBS signal. The constant and PRBS signals were always assigned a norm equal to $\varrho$. 

\subsection{Optimal input design for atomic sensor model}
\label{sec:example_02}
In \cite{JimenezMatinez2018} and also in \cite{dilcher2025}, a simplified paradigmatic model of an atomic sensor (an atomic magnetometer, \cite{FabricantNovikovaBison2023}, \cite{BudkerJacksonKimball2013}) is introduced, in which the dynamics is governed by oscillations of the collective spin of an atomic ensemble subjected to an external magnetic field. The system is driven by circularly polarized light from a pump laser, whose frequency acts as the input signal. A linearly polarized probe laser illuminates the atoms, and upon transmission through the medium, its polarization undergoes a Faraday rotation. The $J_z$ component of the collective spin is inferred from the measurement of the probe laser's polarization angle. The model presented in \cite{JimenezMatinez2018} describes the dynamics of the spin components $\vec{J}=[J_y, J_z]^{\TT}$ and has the form 
\begin{equation}
    d\vec{J}=\begin{bmatrix}
        -\frac{1}{T_2}&\omega_L\\
        -\omega_L&-\frac{1}{T_2}
    \end{bmatrix}Jdt+
    \begin{bmatrix}
        0\\1
    \end{bmatrix}\mathcal{E}(t)dt+d\vec{w}^{(J)},
    \label{eq:equation_1_jimenez}
\end{equation}
where $\omega_L$ is the Larmor frequency, $T_2=0.87$ ms, is the relaxation time, $\mathcal{E}$ is the pumping laser frequency, and $\vec{w}^{(J)}$ is the Wiener process with known covariance $q\mathbb{I}$. The observation has the form $I_k=g_DJ_z(k\Delta)+\xi_k$, $k=0,1, ...$, where $I_k$ is the photocurrent, $\Delta=5$ $\mu\mrm{s}$, is the sampling time, $\xi_k\sim\cN(0,\sigma_{\xi}^2)$ and $g_D$, $\sigma_{\xi}$ are known parameters. The Larmor frequency and the external magnetic field $B$ are related to each other by the formula $\omega_L=\gamma_eB$, where $\gamma_e$ is the gyromagnetic ratio. Hence, by measuring $\omega_L$, one can determine the field $B$. Taking $T_2$ as the time unit and after rescaling the time, state variables, observations, and the input signal $\mathcal{E}$, we get the following, equivalent to \eqref{eq:equation_1_jimenez}, stochastic system:
\begin{equation}
    d\vec{x}=\left(\mat{A}_c(\theta)+\mat{B}_cu\right)dt+\mat{G}_cd\vec{w},
    \label{eq:ct_oscillator}
\end{equation}
\begin{equation}
    y_k=\mat{C}\vec{x}_k+s_vv_k,
\end{equation}
where
\begin{equation}
    \mat{A}_c(\theta)=\begin{bmatrix}
        -1&\theta\\
        -\theta&-1
    \end{bmatrix}, 
    \mat{B}_c=\begin{bmatrix}
        0\\b_c
    \end{bmatrix},
    \mat{G}_c=\sqrt{2}\mathbb{I}, 
    \mat{C}=\begin{bmatrix}
        0&1
    \end{bmatrix},
    \label{eq:sys_matrices_osc_ct}
\end{equation}
$\vec{x}_k=\vec{x}(t_k)$, $t_k=k\Delta$, $\Delta=5.7471\cdot 10^{-3}$, $s_v=11.85$, $b_c=10^5$. The input signal $u$ in \eqref{eq:ct_oscillator} corresponds to $\mathcal{E}$ in \eqref{eq:equation_1_jimenez}. The parameter $\theta$ in \eqref{eq:ct_oscillator} is related to the Larmor frequency $\omega_L$ in \eqref{eq:equation_1_jimenez}, by formula $\theta=\omega_LT_2$. Since the estimation error of the parameter $\theta$ depends on the input signal $u$, a natural question arises as to what form this signal should take. To solve this problem, we will go to discrete time and apply the methodology described in Sections \ref{sec:quasi_linear_sys} and \ref{sec:classical_methods}. Assuming that $u(t)=u_k, t\in [t_k, t_{k+1}]$, the discrete-time system corresponding to \eqref{eq:ct_oscillator} has the form
\begin{equation}
    \vec{x}_{k+1}=\mat{A}(\theta)\vec{x}_k+\mat{B}(\theta)u_k+\mat{G}w_k,
\end{equation}
where, according to the procedure given in Appendix \ref{app:appendix_C}:
\begin{align}
\label{eq:sys_matrices_osc_dt_0}
    &\mat{A}(\theta)=e^{-\Delta}\begin{bmatrix}
        \cos\theta\Delta & \sin\theta\Delta\\
        -\sin\theta\Delta & \cos\theta\Delta
    \end{bmatrix}, 
    \mat{B}(\theta)=\frac{b_c}{1+\theta^2}\begin{bmatrix}
        \theta-e^{-\Delta}(\theta\cos\theta\Delta+\sin\theta\Delta)\\
        1-e^{-\Delta}(\cos\theta\Delta-\theta\sin\theta\Delta)
    \end{bmatrix},\\
    &\mat{G}=\sqrt{1-e^{-2\Delta}}\mathbb{I}.
    \label{eq:sys_matrices_osc_dt}
\end{align}
We assume that the prior distribution of $\theta$ is Gaussian, that is, $p_0(\theta)=\cN(\theta,m_{\theta}, s_{\theta})$ with $m_{\theta}=54.6637$, $s_{\theta}=10.76$, which corresponds to the Larmor frequency of 10 kHz and its initial uncertainty of the order of 600 Hz ($3\sigma$). At the beginning of the process, the system is in thermal equilibrium, that is, $p(\vec{x}_0|\vec{\theta})=\cN(\vec{x}_0, 0, \mathbb{I})$. Following Lemma \ref{lem:Lemma_2}, we get $\theta_1=m_{\theta}-\sqrt{s_{\theta}}\approx 51$, $\theta_2=m_{\theta}+\sqrt{s_{\theta}}\approx 58$, $p_{0,1}=p_{0,2}=0.5$. Thus, $r=2$ in \eqref{eq:inf_lb_formulas_0} and, accordingly to \eqref{eq:eil_two_poss}, the Bayesian optimal signal is a solution of the simplified and convex optimization problem \eqref{eq:opt_problem_two_poss} with $d_{1,2}$ defined by Lemmas \ref{lem:Lemma_3} and \ref{lem:lemma_4}. Moreover, since the matrices in \eqref{eq:sys_matrices_osc_dt_0}, \eqref{eq:sys_matrices_osc_dt} do not depend on $u_k$, the last two terms in \eqref{eq:d_ij_calc} can be omitted. The set of admissible signals is given by \eqref{eq:U_ad} with $\tilde{\vec{U}}=0$, that is, the signal norm cannot be greater than $\varrho$.

In order to employ the classical methods described in Sect. \ref{sec:classical_methods}, it is necessary to first evaluate the sensitivity of the prediction error. Similarly as in the previous example, the polynomials $A, B, C$ have the form 
\begin{align}
 &A(\theta,z)=1-2e^{-\Delta}\cos(\theta\Delta)z^{-1}+e^{-2\Delta}z^{-2},\\ 
 &B(\theta,z)=\mat{B}_2(\theta)-e^{-\Delta}\left(\mat{B}_1(\theta)\sin(\theta\Delta)+\mat{B}_2(\theta)\cos(\theta\Delta)\right)z^{-1},\\
 \begin{split}
 &C(\theta,z)=1+\left(\mat{K}_2(\theta)-2e^{-\Delta}\cos(\theta\Delta)\right)z^{-1}+\\
 &+e^{-\Delta}\left(e^{-\Delta}-\mat{K}_1(\theta)\sin(\theta\Delta)-\mat{K}_2(\theta)\cos(\theta\Delta)\right)z^{-2},
 \end{split}
\label{eq:oscillator_polynomials} 
\end{align}
and the Kalman gain $\mat{K}$ and the vector $\mat{B}$ are given by \eqref{eq:kalman_gain_gh_filters} and \eqref{eq:sys_matrices_osc_dt_0}, respectively. The sensitivity equation \eqref{eq:pem_sensitivity} is given by \eqref{eq:ex_02_sensitivity}. The D-optimal signal is then obtained by maximization of the averaged D-optimal criterion \eqref{eq:av_d_optimal_ex_02} with constraints \eqref{eq:U_ad}.

\begin{figure}[H]
\centering
%\isPreprints{\centering}{} % Only used for preprints
\includegraphics[width=13 cm]{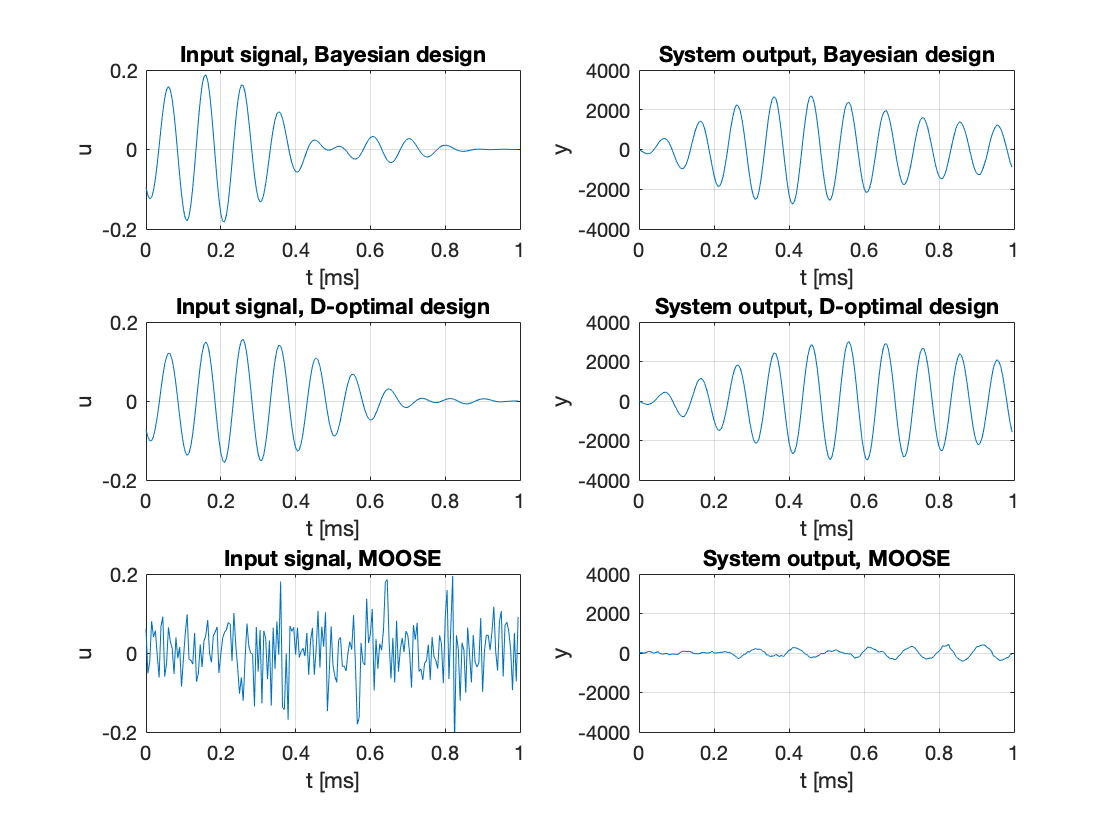}
\caption{Optimal input signals (left) and corresponding system outputs (right) obtained by maximizing the Bayesian criterion \eqref{eq:opt_problem_two_poss} (top), the averaged D-optimal criterion \eqref{eq:av_d_optimal_ex_02} (middle), and the spectral criterion \eqref{eq:pem_fisher_spectral} (bottom), subject to the constraint \eqref{eq:U_ad} with $\tilde{\vec{U}}=0$.
Maximization of the spectral criterion \eqref{eq:pem_fisher_spectral} was performed using the MOOSE-2 solver evaluated at $\theta=m_{\theta}$. The norm of all signals is equal to 1, and the scale is consistent across all plots.
\label{fig:example_02_f1}}
\end{figure}

\begin{figure}[H]
\centering
%\isPreprints{\centering}{} % Only used for preprints
\includegraphics[width=12 cm]{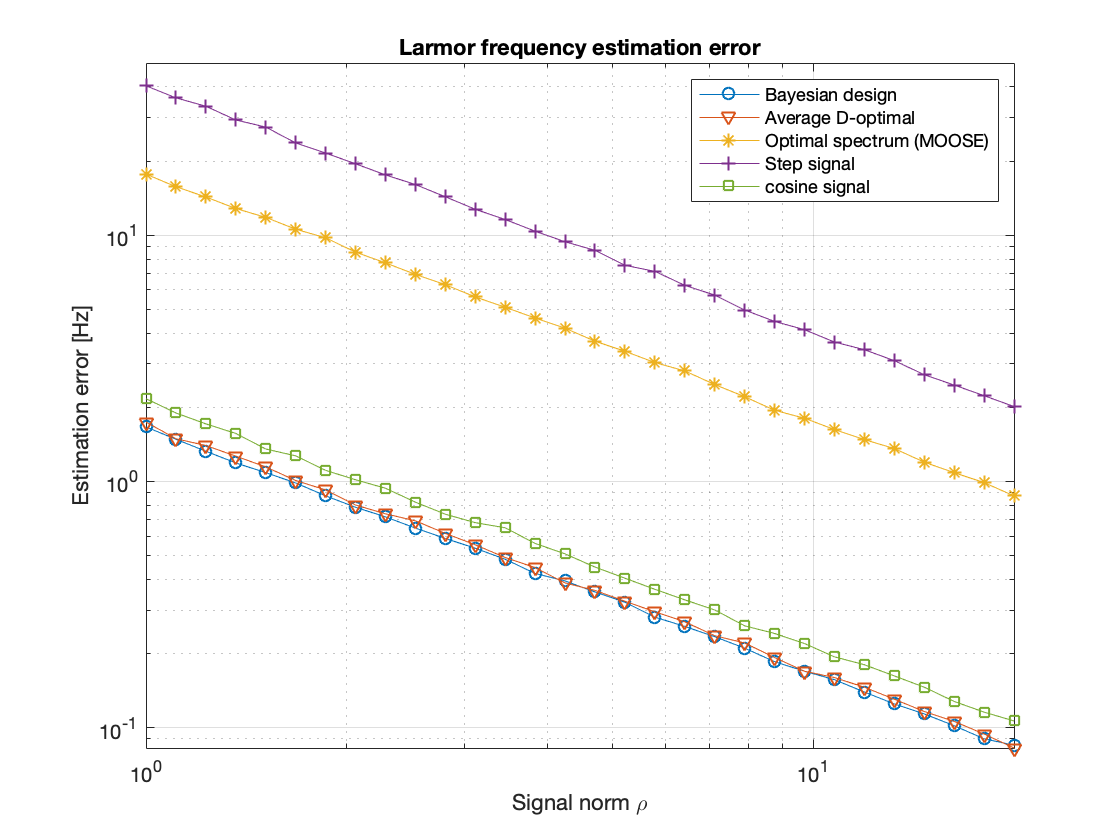}
\caption{Mean estimation errors of the Larmor frequency $f_L=\frac{\omega_L}{2\pi}$, obtained using the MAP estimator \eqref{eq:map_estimator}, as functions of the maximum admissible signal norm $\varrho$, for different signals. The results are based on a Monte Carlo simulation with 3000 repetitions.
\label{fig:example_02_f2}}
\end{figure}

The optimal input signals were designed by maximizing the Bayesian criterion \eqref{eq:inf_lb_simple}, the averaged D-optimal criterion \eqref{eq:av_d_optimal_ex_02}, and the spectral criterion \eqref{eq:pem_fisher_spectral}, subject to the constraint \eqref{eq:U_ad} with $\tilde{\vec{U}}=0$. The results are presented in Figures \ref{fig:example_02_f1} and \ref{fig:example_02_f2}. Figure \ref{fig:example_02_f2} also shows the estimation error for the step (constant) signal and harmonic signal $u(t)=a\cos(m_{\theta}t)$. Frequency of the harmonic signal was equal to the expected value of the a priori distribution of the parameter $\theta$. The constant, MOOSE and harmonic signals were always assigned a norm equal to $\varrho$. 

\subsection{Bayesian input signal design for pump laser control in optically pumped magnetometer}
\label{sec:example_03}
Optically pumped magnetometers operate by aligning atomic spins with a circularly polarized pump laser, after which the spins precess around the external magnetic field at the Larmor frequency. The probe laser measures this precession via polarization rotation (Faraday effect), linking the detected signal to the magnetic field \cite{FabricantNovikovaBison2023}, \cite{BudkerJacksonKimball2013}. The pump laser's frequency strongly affect spin polarization and coherence time making precise laser control central to minimizing estimation error. The advanced control strategies can then suppress the noise and enhance sensitivity. Consequently, accurate control of the pumping laser is a key factor in achieving high-resolution and low-error magnetometer. We consider here the magnetometer model given by equation (S9) in the article \cite{Truullinou2021}: 
\begin{equation}
\frac{d\vec{F}}{d\tau}=\left(-\gamma_e\vec{B}+GS_3\hat{z}\right)\times \vec{F}-\gamma \vec{F}+P(\tau)(\hat{z}F_{max}-\vec{F})+\mat{G}_0(P(\tau))\vec{w},
\label{eq:model_S9}
\end{equation}
where $\vec{F}=(F_x, F_y, F_z)^{\TT}$ is the collective atomic spin, $\gamma_e$ is electron gyromagnetic ratio, $\vec{B}=(B_x, B_y, B_z)^{\TT}$ is constant magnetic field vector, $G$ is known positive constant, $GS_3\hat{z}$ is the effective field produced by ac-Stark shifts due to the probe laser where $S_3$ is white Gaussian noise with variance $\sigma_3^2$. The optical pumping rate $P(\tau)\geqslant 0$ is an input signal. The atomic spin noise $\mat{G}_0(P(\tau))\vec{w}$ is modelled as white Gaussian where $\vec{w}=(w_x, w_y, w_z)^{\TT}$ is a vector of standard and mutually independent Wiener increments. The $\mat{G}_0$ matrix is diagonal and is given by
\begin{equation}
\mat{G}_0(P(\tau))=\sqrt{\tfrac{2}{3}F(F+1)N_A(\gamma+P(\tau))}\mathbb{I},
\label{eq:G0_matrix_S9}
\end{equation}
where $N_A$ is the number of atoms and  $F$ is known atomic spin number. Parameter $F_{max}=N_AF$ is the maximum possible polarization. The transverse relaxation rate $\gamma$ depends on the number of atoms and is given by $\gamma(N_A)=T_2(N_A)^{-1}=\gamma_0+10^{-12}\alpha N_A$, where $\gamma_0$ and $\alpha$ are known positive constants and $T_2$ is the effective coherence time. The observation equation has the form
\begin{equation}
S_2=F_z+N_{S_2},
\label{eq:observation_S9}
\end{equation}
where $N_{S_2}$ denote measurement noise with variance $\sigma_2^2$. In the experiment, the $S_2$ component of the Stokes vector is measured at discrete time moments $t_k=k\Delta$, where $\Delta$ is sampling period. The realistic parameters of the model are given in Table \ref{tab:table_1}.

\begin{table}[H]
\centering
\caption{Typical parameters.}
\label{tab:table_1}

\begin{tabular}{|c|c|c|}
\hline
\textbf{Parameter} & \textbf{Abbreviation} & \textbf{Typical value} \\ 
\hline
Number of atoms & $N_A$ & $10^{12}$ \\ 
\hline
Spin number & $F$ & $1$ \\
\hline 
Larmor frequencies & $\gamma_e\vec{B}$ &  $2\pi [-50,50]$ kHz\\
\hline 
Parameter & $\gamma_0$ & 600 Hz\\
\hline
Parameter & $\alpha$ & 550 Hz\\
\hline
Typical relaxation time & $T_2$ & 0.87 ms \\ 
\hline
Typical relaxation rate & $\gamma$ & 1149 Hz \\
\hline
Pumping rate & $P$ & 0-200 kHz \\
\hline  
Measurement noise level & $\sigma_2$ & $9.6755\cdot 10^{6}$ \\
\hline 
Sampling time & $\Delta$ & $5\mu\mathrm{s}$\\
\hline
\end{tabular}
\end{table}

In what follows, equation \eqref{eq:model_S9} will be interpreted in the It\"o sense. \textbf{Moreover, we assume that the noise $GS_3$ in \eqref{eq:model_S9} is small and can be omitted.} By introducing state variables $\vec{\xi}=\vec{F}\sqrt{\frac{3}{F(F+1)N_A}}$, control variable $u=P/\gamma$ and non-dimensional time $t=\gamma \tau$, and after multiplying both sides of \eqref{eq:observation_S9} by $\sqrt{\frac{3}{F(F+1)N_A}}$, we get the following model:
\begin{equation}
d\vec{\xi}=(\mat{A}_c(\vec{\theta},u)\vec{\xi}+\mat{B}_cu)dt+\mat{G}_c(u)d\vec{\eta}, y_k=\vec{\xi}_3(t_k)+\sigma_v v_k,
\label{eq:ct_model_ex_03}
\end{equation}
where $\vec{\eta}$ is the three-dimensional standard Wiener process with unit covariance and
\begin{equation}
\mat{A}_c(\vec{\theta},u)=\begin{bmatrix}
-(1+u)&\theta_3&-\theta_2\\
-\theta_3&-(1+u)&\theta_1\\
\theta_2&-\theta_1&-(1+u)
\end{bmatrix}, 
\mat{B}_c=\begin{bmatrix}
0\\0\\b_c
\end{bmatrix},
\mat{G}_c(u)=\sqrt{2(1+u)}\mathbb{I},
\label{eq:ABCG_ex_03}
\end{equation}
with $b_c=\sqrt{\tfrac{3N_AF}{F+1}}$, $\sigma_v=\sigma_2\sqrt{\frac{3}{F(F+1)N_A}}$. Taking the parameters from Table \ref{tab:table_1} we have $b_c=1.22\cdot 10^6$, $\sigma_v=11.85$. The parameter vector $\vec{\theta}=(\theta_1, \theta_2, \theta_3)^{\TT}$, represents the external magnetic field due to the relation $\vec{\theta}=\gamma_eT_2\vec{B}$.
If $u$ is a constant signal, then system \eqref{eq:ct_model_ex_03} approaches thermodynamic equilibrium, with $\E \vec{\xi}(t)=-\mat{A}_c(\vec{\theta},u)^{-1}\mat{B}_c u$ and $\mrm{cov}(\vec{\xi}(t))=\mathbb{I}$.

A closer examination of equation \eqref{eq:ct_model_ex_03} shows that the component $\vec{\xi}_3(t,\vec{\theta})$ of its solution remains invariant under rotations of the vector $\vec{\theta}$ about the $z$-axis. As a result, $\vec{\theta}$, and hence the field $\vec{B}$, cannot be uniquely identified from the observations $y_0, \ldots, y_N$. The only quantities that can be uniquely identified in this model are the magnitude of the vector $\vec{B}$ and the angle $\eta$ between $\vec{B}$ and one of the coordinate axes, say the $\hat{z}$ axis. However, to simplify the problem as much as possible, we introduce here the additional assumption that the field $\vec{B}$ always lies in the $x$-$y$ plane, that is, $\vec{B}=(B_x, B_y, 0)^{\TT}$. With this assumption, the change of variables
\begin{equation}
    x_1=\xi_1\sin\varphi-\xi_2\cos\varphi, x_2=\xi_3,
\end{equation}
\begin{equation}
    \cos\varphi=\frac{\theta_1}{\sqrt{\theta_1^2+\theta_2^2}}, \sin\varphi=\frac{\theta_2}{\sqrt{\theta_1^2+\theta_2^2}},
\end{equation}
reduces model \eqref{eq:ct_model_ex_03} to a two-dimensional system:
\begin{equation}
d\vec{x}=(\mat{A}_c(\theta,u)x+\mat{B}_cu)dt+\mat{G}_c(u)d\vec{w}, y_k=\mat{C}\vec{x}_k+\sigma_v v_k,
\label{eq:ct_model_2_ex_03}
\end{equation}
where $\vec{x}=(x_1, x_2)^{\TT}$, $\theta=\gamma_eT_2\sqrt{B_x^2+B_y^2}$, $\vec{w}$ is a two-dimensional standard Wiener process with unit covariance and
\begin{equation}
\mat{A}_c(\theta,u)=\begin{bmatrix}
-(1+u)&\theta\\
-\theta&-(1+u)\\
\end{bmatrix}, 
\mat{B}_c=\begin{bmatrix}
0\\b_c
\end{bmatrix}, 
\mat{C}=\begin{bmatrix}
0&1
\end{bmatrix}, 
\mat{G}_c(u)=\sqrt{2(1+u)}\mathbb{I}.
\label{eq:ABCG_2_ex_03}
\end{equation}
Hence, under the assumption $B_z=0$, the observations $y_0, \ldots, y_N$, variable $\vec{\xi}_3$, and the $\vec{F}_z$ component of the collective spin, are fully characterized by the reduced model \eqref{eq:ct_model_2_ex_03}. Furthermore, within this reduced model it can be readily verified that $\theta$ is uniquely identifiable. Naturally, the accuracy of estimating $\theta$ depends on the choice of input $u$. To determine an input $u$ that maximizes the information about $\theta$, we now turn to the discrete-time formulation of \eqref{eq:ct_model_2_ex_03} and apply the methods described in Sect. \ref{sec:appr_solutions} and \ref{sec:quasi_linear_sys}. Assuming the control signal is piecewise constant, that is, $u(t)=u_k, t\in [t_k, t_{k+1}]$, the process $\vec{x}_k=\vec{x}(t_k)$ satisfies the difference equation
\begin{equation}
\vec{x}_{k+1}=\mat{A}(\theta,u_k)\vec{x}_k+\mat{B}(\theta, u_k)+\mat{G}(u_k)\vec{w}_k,
\label{eq:dt_dynamics_simply}
\end{equation}
where $\vec{w}_k\sim\cN(0, \mathbb{I})$ and the matrices $\mat{A}$, $\mat{B}$, $\mat{G}$ can be calculated following the procedure given in Appendix \ref{app:appendix_C}. Upon completion of straightforward calculations, we get: 
\begin{align}
\label{eq:sys_matrices_2_ex_03}
    &\mat{A}(\theta,u_k)=e^{-(1+u_k)\Delta}\begin{bmatrix}
        \cos\theta\Delta & \sin\theta\Delta\\
        -\sin\theta\Delta & \cos\theta\Delta
    \end{bmatrix},\mat{G}(u_k)=\sqrt{1-e^{-2(1+u_k)\Delta}}\mathbb{I},\\
    &\mat{B}(\theta,u_k)=\frac{b_cu_k}{(1+u_k)^2+\theta^2}\begin{bmatrix}
        \theta-e^{-(1+u_k)\Delta}\left(\theta\cos(\theta\Delta)+(1+u_k)\sin(\theta\Delta)\right)\\
        e^{-(1+u_k)\Delta}\left(\theta\sin(\theta\Delta)-(1+u_k)\cos(\theta\Delta)\right)+(1+u_k)
    \end{bmatrix}.
\end{align}
At the beginning of the process, the system is in a thermal equilibrium corresponding to $u\equiv0$. Hence, $p(\vec{x}_0|\vec{\theta})=\cN(\vec{x}_0, 0, \mathbb{I})$. We also assume that the prior distribution of $\theta$ is Gaussian, that is, $p_0(\theta)=\cN(\theta,m_{\theta}, s_{\theta})$ with $m_{\theta}=54.6637$, $s_{\theta}=3\cdot 10^{-3}$, which corresponds to the Larmor frequency of 10 kHz and its initial uncertainty of the order of 30 Hz ($3\sigma$). Similarly as in Sect. \ref{sec:example_02}, we get from Lemma \ref{lem:Lemma_2}: $\theta_1=m_{\theta}-\sqrt{s_{\theta}}\approx 54.61$, $\theta_2=m_{\theta}+\sqrt{s_{\theta}}\approx 54.72$, $p_{0,1}=p_{0,2}=0.5$. Since $r=2$ in \eqref{eq:inf_lb_formulas_0}, then accordingly to \eqref{eq:eil_two_poss}, the Bayesian optimal signal is a solution of the simplified optimization problem \eqref{eq:opt_problem_two_poss} with $d_{1,2}$ calculated by Lemmas \ref{lem:Lemma_3} and \ref{lem:lemma_4}.
Unlike in the previous examples, in this problem we maximize criterion \eqref{eq:opt_problem_two_poss} with constraints on the signal amplitude, that is, $0\leqslant u_k\leqslant u_{max}$, which is preferable in realistic scenarios.

\begin{figure}
\centering
%\isPreprints{\centering}{} % Only used for preprints
\includegraphics[width=10 cm]{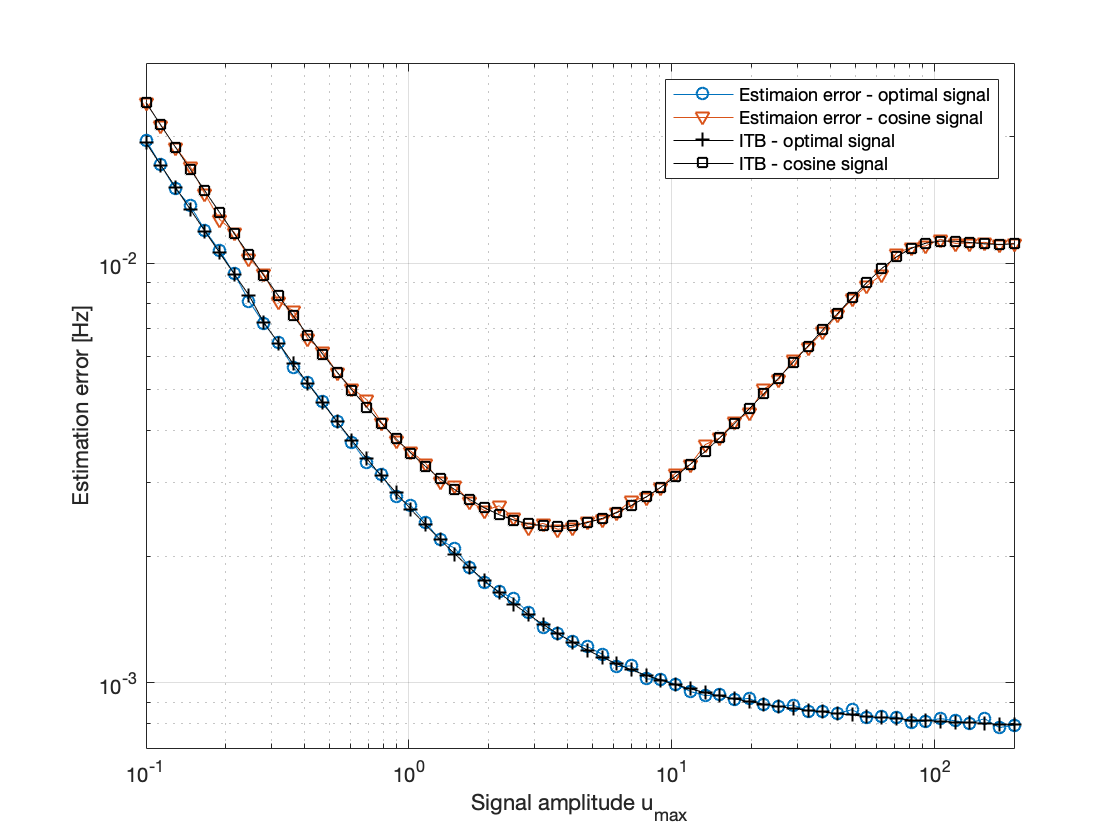}
\caption{The estimation error of the Larmor frequency $f_L = \tfrac{\theta}{2\pi T_2}$ and the Information-Theoretic Bound (ITB), \eqref{eq:error_lb} as a function of the maximum admissible signal amplitude $u_{\max}$. The errors were computed using the MAP estimator \eqref{eq:map_estimator}. Both the errors and the ITB were calculated for two cases: (i) the optimal input signal and (ii) for the the harmonic input $u(t) = 0.5u_{\max}(1+\cos(m_{\theta}t))$. Results are based on a Monte Carlo simulation with 2000 repetitions. The prior was Gaussian with mean Larmor frequency $f_L=10$ kHz, and with its initial uncertainty $\sigma_{f_L}=10$ Hz.
\label{fig:example_03_f1}}
\end{figure}

\begin{figure}
\centering
\includegraphics[width=12 cm]{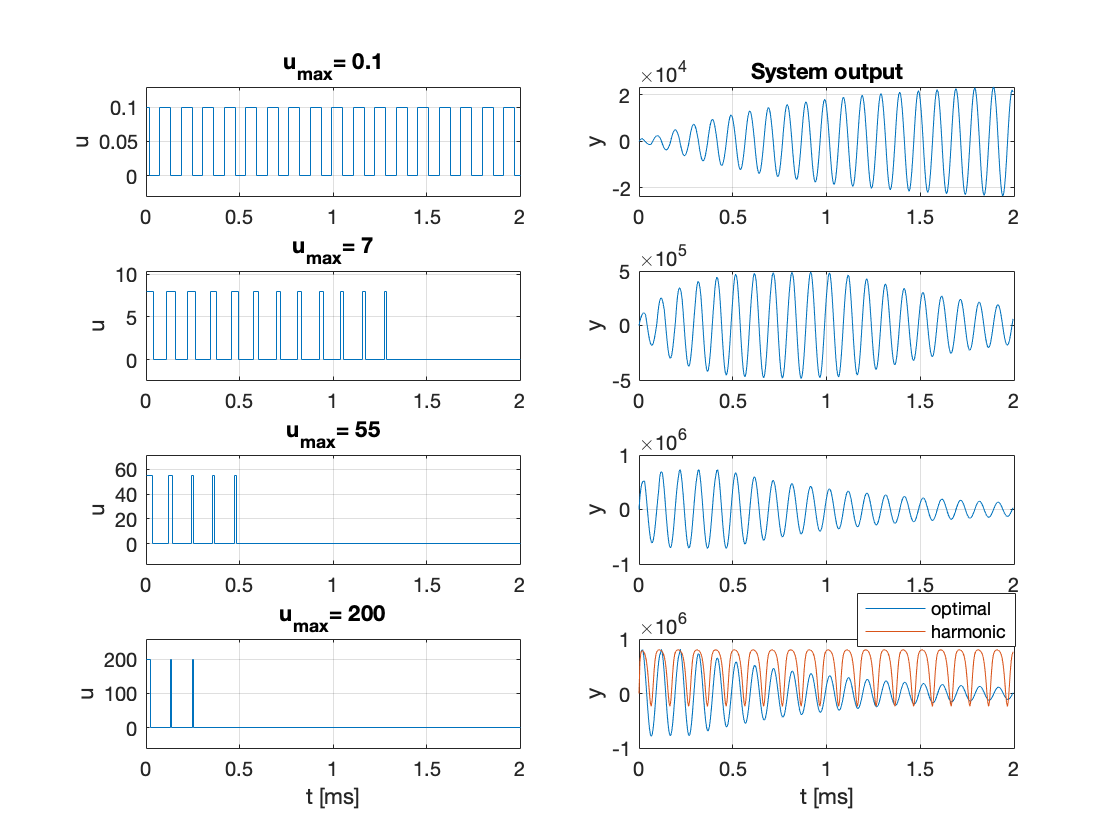}
\caption{The optimal signals with small medium and large amplitude and the corresponding system outputs. The figure in the lower right panel also shows the system output for the harmonic signal $u(t)=0.5u_{max}(1+\cos(m_{\theta}t))$. The prior was Gaussian with mean Larmor frequency $f_L=10$ kHz, and with its initial uncertainty $\sigma_{f_L}=10$ Hz.
\label{fig:example_03_f2}}
\end{figure}

\begin{figure}
\centering
\includegraphics[width=12 cm]{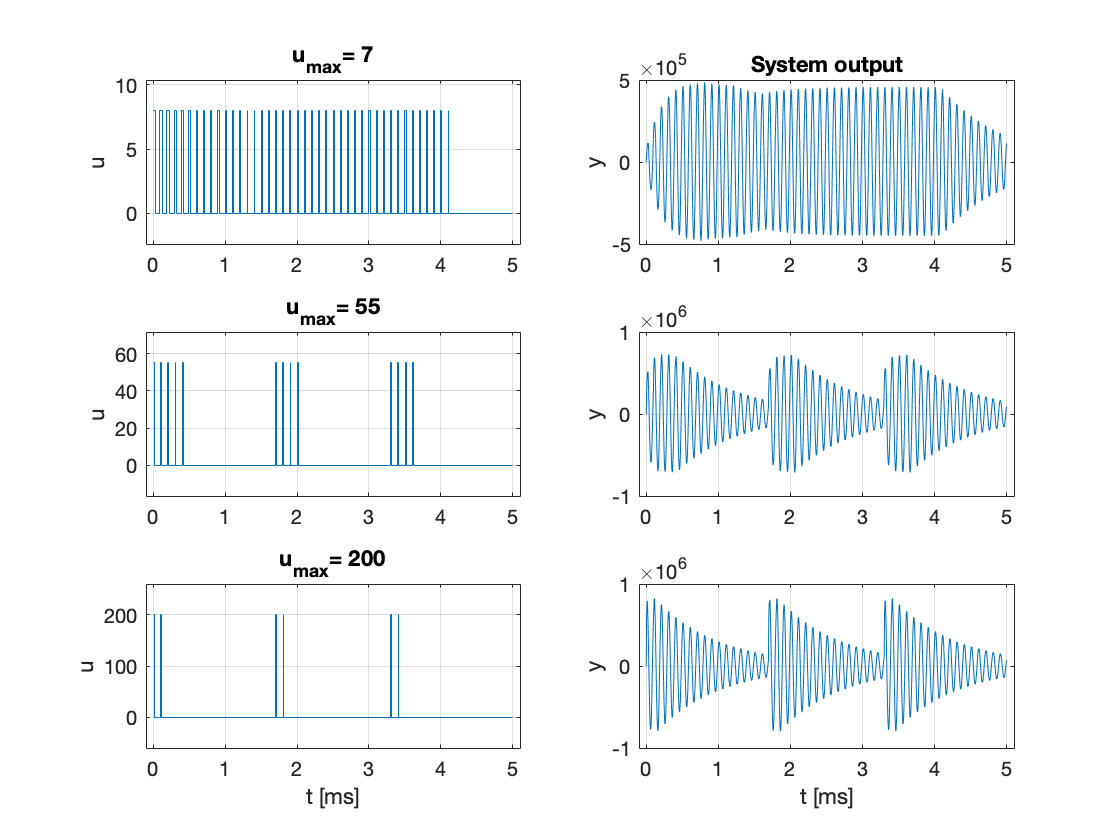}
\caption{Optimal input signals of small, medium, and large amplitudes with corresponding system outputs for a 5 ms experiment.
The prior was Gaussian with mean Larmor frequency $f_L=10$ kHz, and with its initial uncertainty $\sigma_{f_L}=10$ Hz.
\label{fig:example_03_f3}}
\end{figure}

\begin{figure}
\centering
\includegraphics[width=12 cm]{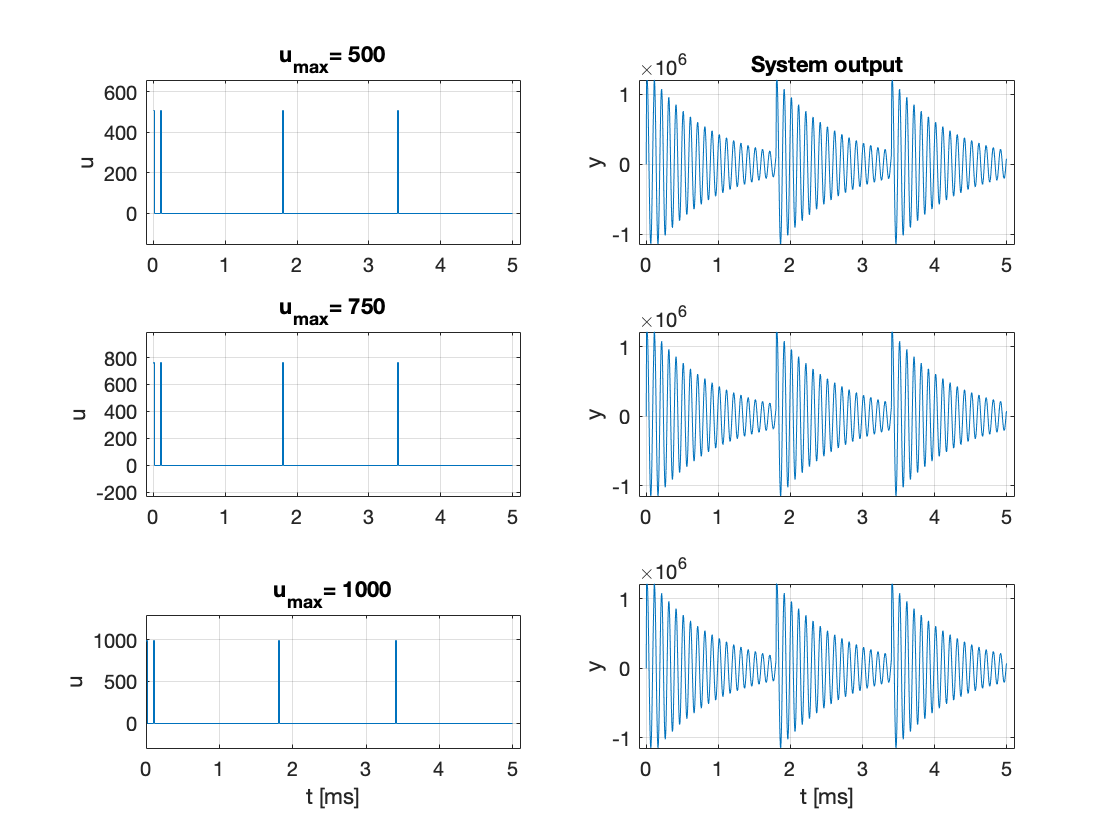}
\caption{Optimal input signals with very large amplitudes and corresponding system outputs. Continuation of Fig. \ref{fig:example_03_f3}.
\label{fig:example_03_f4}}
\end{figure}

The results are presented in Figures~\ref{fig:example_03_f1}-\ref{fig:example_03_f4}. 
The optimal input signal consistently lies on the boundary of the admissible set. 
For small $u_{\max}$, the optimal signal is rectangular, with a frequency close to the \textit{a priori} Larmor frequency. 
Since large values of $u(t)$ strongly damp spin oscillations and increases the noise, the optimal signals for large $u_{\max}$ consist of short pulses at the maximum admissible amplitude. 
Once the oscillations decay, the system should be re-excited by a new sequence of short pulses, repeated periodically, as illustrated in Fig.~\ref{fig:example_03_f3}. 
The harmonic signal $u(t) = 0.5u_{\max}(1+\cos(m_{\theta}t))$ is nearly optimal for small $u_{\max}$ but becomes ineffective for large $u_{\max}$, as it strongly damps the oscillations (see the lower right panel of Fig.~\ref{fig:example_03_f2}). 
As a result, the measurements carry less information about the Larmor frequency, and the estimation error increases despite the higher signal amplitude. An analogous behaviour is observed for rectangular signals. More generally, let $s(t)\in [0,1]$, be any signal and define $u(t)=as(t)$ with $a\geqslant 0$. Then, as illustrated in Fig.~\ref{fig:example_03_f1}, estimation error reaches a minimum for some non-zero value of the parameter $a$.

Extending the experimental duration from $2$ to $5$~ms reduces the estimation error by a factor of 2 compared to the case shown in Fig.~\ref{fig:example_03_f1}. For $u_{\max}=200$ and an experiment duration of $5$ ms, the harmonic input signal yields an estimation error of $7$ mHz, while the optimal signal, shown in the lower left panel of Fig.~\ref{fig:example_03_f3}, reduces the error to $0.48$ mHz, that is, approximately $14$ times smaller. Further increasing the maximum allowable amplitude leads to a signal that consists of $5\mu\rm{s}$ pulses of maximum amplitude, spaced approximately $1.5T_2$ apart. This is illustrated in Fig. ~\ref{fig:example_03_f4}.

Finally, the estimation error attains the Information-Theoretic Lower Bound~\eqref{eq:error_lb}, demonstrating that in this case the MAP estimator~\eqref{eq:map_estimator} achieves optimal performance.

It should be noted that the above models assume the Markovian environment and this condition should be checked in an experiment. To do this, one can use the criterion given in \cite{Breuer2009}. Non-Markovian models are much more complicated (see e.g. \cite{Shen2018}) and one would need to employ a noise model with long memory. To model long-memory noise, fractional-order stochastic equations can be used instead of \eqref{eq:model_S9}. Such models capture long-memory effects, and their noise correlation function decays slowly, for example as $t^{-1/2}$, \cite{Iomin_2009}, \cite{Bania2013_Laguerre}, \cite{Bania_Laguerre}.

To implement the proposed method in real time, one can proceed as follows. First, observe that the pump signal has a simple structure, consisting of short pulses at the maximum admissible amplitude, each lasting approximately \(5\,\mu\text{s}\). These pulses should be repeated with a period of roughly \(1.5T_2\), and each pulse should be triggered when the vector \([F_y, F_z]\) forms an angle of about \(30^\circ\) with the \(z\)-axis (i.e., \(30^\circ\) before the maximum of \(F_z\)). To estimate the unknown vector \(\vec{F}\) and the Larmor frequency, the MAP estimator is too slow for real-time applications. Instead, an Extended Kalman Filter (EKF) can be employed in a manner roughly similar to that described in  \cite{dilcher2025},  \cite{Magrini2021} and \cite{Amoros2025}. This approach is considered feasible for implementation in an experimental setup.


%*******************
\section{Discussion and conclusions}
This paper has developed a Bayesian framework for optimal input signal design in the identification of quasi-linear stochastic dynamical systems. By an Information-Theoretic Lower Bound on estimation error and its connection to the Bayesian Cram\'{e}r-Rao Bound, we showed that maximizing mutual information provides a principled alternative to Fisher-information-based criteria. The proposed method relies on the maximization of the MI lower bound \eqref{eq:inf_lb}, which produces a tractable surrogate objective for both finite parameter sets and parameter spaces of continuum cardinality. A key contribution is the algorithmic reduction of the dimension of covariance matrices required for inversion by a factor of $N$, making the method feasible for long-term experiments.

The comparison with the average D-optimal design highlights the practical benefits of the Bayesian approach. While classical methods are computationally efficient, they require complex differentiations to evaluate parameter sensitivities and may yield suboptimal results when parameter uncertainty is large or when the system exhibits significant non-linearities. In contrast, the proposed Bayesian method requires only the system matrices $\mat{A}$, $\mat{B}$, $\mat{C}$, $\mat{G}$, 
together with the prior distributions of the parameter and initial conditions, without the need to calculate derivatives of prediction errors. This makes the method applicable to a much broader class of systems, while also enabling it to handle large initial parameter uncertainty.

The method also has certain limitations. The lower bound of the MI involves exponential terms that can vanish when the pairwise distance factors $d_{i,j}(\vec{U})$ are large, which can cause numerical problems. However, this drawback can be mitigated by appropriate scaling of the optimization problem. If we consider the simplified optimization problem \eqref{eq:opt_problem_two_poss}, with only two candidate parameter values, these numerical problems never occur. 
The discrete approximation of the MI \eqref{eq:MI_discrete} is a potential source of problems and the weights and nodes in \eqref{eq:finite_gauss_mix_appr} should be carefully selected to achieve sufficient approximation accuracy.
The third limitation arises from the fact that the maximized criterion is only a lower bound on the MI and is generally not tight. Consequently, there certainly exists a class of problems for which maximizing this lower bound is inefficient and may generate signals that are far from optimality in the sense of maximizing the MI \eqref{eq:MI}.

In all analyzed examples, the proposed Bayesian approach, although approximate, generated signals no worse, or even better (see Figs. \ref{fig:example_01_f2} and \ref{fig:example_02a_f2}), than the classical methods. 

The second example illustrates that a non-Gaussian prior distribution leads to increased errors in the average D-optimal method. For a Gaussian prior distribution, it was confirmed that both the average D-optimal and the proposed Bayesian method yield identical results. This observation underscores the sensitivity of the classical approach to the form of the prior distribution and highlights the necessity of employing estimation techniques that are robust to non-Gaussian priors. In the third example, the D-optimal method produces results almost identical to those of the Bayesian approach. To explain this, note that in this problem the prior distribution of the parameter $\theta$ is relatively narrow. Then the function $d_{1,2}(\vec{U})$, which we minimize in this task, is approximately proportional to the sensitivity of the output to changes in $\theta$. Thus, $d_{1,2}(\vec{U})$ can be interpreted as a quantity proportional to the Fisher information. Consequently, the resulting input signals and the corresponding estimation errors are nearly identical.

The study of atomic sensor models further demonstrates the practical relevance of the approach. The optimal signals in these examples are always better than the harmonic signal with a frequency equal to the expected natural frequency of the oscillator. The fourth example, a seemingly minor modification of the oscillator from the third example, shows that the dependence of the system matrices on the control signal is significant and leads to completely different optimal signals.
In the analyzed examples, the MAP estimator achieves an Information-Theoretic Lower Bound \eqref{eq:error_lb}, but this is not always the case and, depending on the task, there are better estimators. Unfortunately, finding them is difficult.

Since the method produces the posterior distribution $p(\vec{\theta}, \vec{x}_k|\vec{Y}_k)$, it can be easily converted to a sequential Bayesian Adaptive Design (BAD) algorithm \cite{RainforthEtAl2024}, \cite{HuanJagalurMarzouk2024}. Then the optimal strategy is a functional of the posterior, that is, $\vec{u}_k=\vec{\phi}_k(p(\vec{\theta}|\vec{Y}_k), \vec{m}_k(\vec{\theta}), \vec{S}_k(\vec{\theta}))$. In the simplest case, when the matrices $\mat{A}$, $\mat{B}$, $\mat{G}$ do not depend on $\vec{u}_k$ and $\vec{\Theta}=\lbrace\vec{\theta_1}, \vec{\theta}_2\rbrace$, the optimal strategy $\vec{\phi}_k$ can be determined by maximizing \eqref{eq:inf_lb} on the trajectories of the system \eqref{eq:kalman_ij_0}-\eqref{eq:kalman_ij_1}. This problem is deterministic, therefore, relatively simple, and can be solved using deterministic optimal control methods.

From a broader perspective, quasi-linear systems arise naturally in quantum mechanics, chemical engineering, and thermal processes, making the proposed method widely applicable.  In conclusion, this work provides both theoretical justification and practical tools for Bayesian input design in quasi-linear stochastic systems. By bridging information-theoretic principles with efficient computational methods, it establishes a foundation for robust experimental design in a wide range of applications. The results reported here should stimulate further research at the intersection of Bayesian inference, control, and the identification of non-linear systems.

\bibliography{references_v5.bib}

\begin{appendices}

\section{Proofs}
\label{app:appendix_A}
\begin{proof}[Proof of Theorem 1.]
Let $\hat{\vec{\theta}}_{M}(\vec{Y},\vec{U})=\E_{p(\vec{\theta}|\vec{Y},\vec{U})}(\vec{\theta}|\vec{Y},\vec{U})$ be the Minimum Mean Squared Error estimator of $\vec{\theta}$.  
The conditional covariance of $\hat{\vec{\theta}}_M$ is given by $\mat{C}(\vec{Y},\vec{U})=\int p(\vec{\theta}|\vec{Y},\vec{U})(\vec{\theta}-\hat{\vec{\theta}}_M(\vec{Y},\vec{U}))(\vec{\theta}-\hat{\vec{\theta}}_M(\vec{Y},\vec{U}))^{\TT}d\vec{\theta}$. Since Gaussian distribution maximizes entropy over all distributions with the same covariance, it can be proved, (see \cite{CoverThomas2006}[Thm. 8.6.5,  p. 255]) that   
\begin{equation}
H_{\theta|Y}(\vec{U})=\E(-\ln p(\vec{\theta}|\vec{Y},\vec{U}))\leqslant \tfrac{1}{2} \E\ln ((2\pi e)^{n_\theta}|\mat{C}(\vec{Y},\vec{U})|).
\label{eq:lb_proof_01}
\end{equation}
Any covariance matrix $\mat{C}$ satisfies the inequality $|\mat{C}|\leqslant (n_{\theta}^{-1}\mathrm{tr}(\mat{C}))^{n_{\theta}}$ (see \cite{CoverThomas2006}, Thm. 17.9.4, p. 680). Hence  
\begin{equation}
\ln |\mat{C}(\vec{Y},\vec{U})|\leqslant n_{\theta}\ln n_{\theta}^{-1}\mathrm{tr}\mat{C}(\vec{Y},\vec{U}).
\label{eq:lb_proof_02}
\end{equation}
Taking into account \eqref{eq:MI} and \eqnsref{eq:lb_proof_01}{eq:lb_proof_02} we have
\begin{align*}
& H_{\theta}-I_{\theta;Y}(\vec{U})=H_{\theta|Y}(\vec{U})=\E (-\ln p(\theta |\vec{Y},\vec{U})) \leqslant \\
& \tfrac{1}{2} \E\ln ((2\pi e)^{n_\theta}|\mat{C}(\vec{Y},\vec{U})|)\leqslant \tfrac{n_{\theta}}{2}\E\ln (2\pi e n_{\theta}^{-1}\mathrm{tr}\mat{C}(\vec{Y},\vec{U})).
\end{align*}
By concavity of the logarithm and from Jensen's inequality
\begin{equation*}
H_{\theta}-I_{\theta;Y}(\vec{U})\leqslant\tfrac{n_{\theta}}{2}\ln (2\pi e n_{\theta}^{-1}\E\mathrm{tr}\mat{C}(\vec{Y},\vec{U})).
\end{equation*}
Using the equality $\E\mathrm{tr}\mat{C}(\vec{Y},\vec{U})=\E |\vec{\theta} - \hat{\vec{\theta}}_{M}(\vec{Y},\vec{U}) |^2$,
yields
\begin{equation*}
H_{\theta}-I_{\theta;Y}(\vec{U})\leqslant\tfrac{n_{\theta}}{2}\ln (2\pi e n_{\theta}^{-1}\E |\vec{\theta} - \hat{\vec{\theta}}_{M}(\vec{Y},\vec{U}) |^2).
\end{equation*}
Since $\hat{\vec{\theta}}_M$ is MMSE then
$\E |\vec{\theta} - \hat{\vec{\theta}}_{M}(\vec{Y},\vec{U}) |^2\leqslant \E |\vec{\theta} - \hat{\vec{\theta}}(\vec{Y},\vec{U}) |^2$ and 
\begin{equation*}
H_{\theta}-I_{\theta;Y}(\vec{U})\leqslant\tfrac{n_{\theta}}{2}\ln (2\pi e n_{\theta}^{-1}\E |\vec{\theta} - \hat{\vec{\theta}}(\vec{Y},\vec{U}) |^2),
\end{equation*}
which is equivalent to the first inequality in \eqref{eq:error_lb}. The proof of the Efroimovitch inequality, that is, the second inequality in \eqref{eq:error_lb}, is given in \cite{Lee2022}[Cor. 3, Ch. 2.2, p. 16].
\end{proof}
%%%% Lemma 3 proof ***************************************************************************************************
\begin{proof}[Proof of Lemma 3]
The proof of \eqref{eq:NYF_calc}, \eqref{eq:pem_kalman_0}-\eqref{eq:pem_kalman} is well documented in the literature and can be found in \cite{BanBar2016} and \cite{Sarkka2013}[Thm. 12.3, p. 187]. However, for completeness and the convenience of the reader, we reproduce it here in its entirety. Let us denote $\vec{X}_k=\mrm{col}(\vec{x}_0, \dots, \vec{x}_k)$, $\vec{Y}_k=\mrm{col}(\vec{x}_0, \dots, \vec{x}_k)$. We begin by recalling the filtering equations. If $\vec{\theta}$ is fixed parameter, then the solution of equation \eqref{eq:sys_dynamics} is a Gauss-Markov process with transition density
\begin{equation}
p(\vec{x}_k|\vec{x}_{k-1},\vec{\theta})=\cN(\vec{x}_k,\mat{A}_{k-1}\vec{x}_{k-1}+\mat{B}_{k-1}, \mat{A}_{k-1}\mat{S}_{k-1}\mat{A}_{k-1}^{\TT}+\mat{G}_{k-1}\mat{G}_{k-1}^{\TT}),
\label{eq:lem3p_transition_pdf}
\end{equation}
where we used the notation of Sect.\ref{sec:quasi_linear_sys} and we will omit the arguments $\vec{U}$ and $\vec{u}_k$, in all the formulas below. It follows from \eqref{eq:sys_observation}, that the density of the observations $\vec{y}_k$, conditioned on $\vec{X}_k$, $\vec{Y}_{k-1}$, $\vec{\theta}$, has the form:   
\begin{equation}
p(\vec{y}_k|\vec{X}_k,\vec{Y}_{k-1},\vec{\theta})=p(\vec{y}_k|\vec{x}_k, \vec{\theta})=\cN(\vec{y}_k,\mat{C}\vec{x}_k,\mat{S}_v(\vec{\theta})).
\label{eq:lem3p_y_pdf}
\end{equation}
According to the assumptions given at the beginning of Sect. \ref{sec:quasi_linear_sys}, the initial distribution of $\vec{x}_0$ is given by: 
\begin{equation}
p(\vec{x}_0|\vec{\theta})=\cN(\vec{x}_0,\vec{m}_0^-(\vec{\theta}),\mat{S}_0^-(\vec{\theta})),
\label{eq:lem3p_xini_pdf}
\end{equation}
where $\vec{m}_0^-$, $\mat{S}_0^-$ are smooth functions and $\mat{S}_0^-(\vec{\theta})>0$, for all $\vec{\theta}\in\vec{\Theta}$. To find $p(\vec{Y}_k|\vec{\theta})$ we proceed as follows:
\begin{equation}
\begin{split}
&p(\vec{x}_k, \vec{Y}_k|\vec{\theta})=\int p(\vec{y}_k,\vec{x}_k, \vec{x}_{k-1}, \vec{Y}_{k-1}|\vec{\theta})d\vec{x}_{k-1}=\\
&=\int p(\vec{y}_k|\vec{x}_k, \vec{x}_{k-1}, \vec{Y}_{k-1}, \vec{\theta})p(\vec{x}_k, \vec{x}_{k-1}, \vec{Y}_{k-1}| \vec{\theta})  d\vec{x}_{k-1}=\\
&=p(\vec{y}_k|\vec{x}_k, \vec{\theta})\int p(\vec{x}_k| \vec{x}_{k-1}, \vec{Y}_{k-1}, \vec{\theta})p(\vec{x}_{k-1},\vec{Y}_{k-1}, \vec{\theta})d\vec{x}_{k-1}=\\
&=p(\vec{y}_k|\vec{x}_k, \vec{\theta})\int p(\vec{x}_k| \vec{x}_{k-1}, \vec{\theta})p(\vec{x}_{k-1}|\vec{Y}_{k-1}, \vec{\theta})d\vec{x}_{k-1}p(\vec{Y}_{k-1}|\vec{\theta})=\\
&=p(\vec{y}_k|\vec{x}_k, \vec{\theta})p(\vec{x}_k| \vec{Y}_{k-1}, \vec{\theta})p(\vec{Y}_{k-1}|\vec{\theta}),
\end{split}
\label{eq:lem3p_filter_der}
\end{equation}
where
\begin{equation}
p(\vec{x}_k| \vec{Y}_{k-1}, \vec{\theta})=\int p(\vec{x}_k| \vec{x}_{k-1}, \vec{\theta})p(\vec{x}_{k-1}|\vec{Y}_{k-1}, \vec{\theta})d\vec{x}_{k-1},
\label{eq:lem3p_pred_distr}
\end{equation}
is the so called predictive distribution or prediction step. Integration of both sides of \eqref{eq:lem3p_filter_der} over $\vec{x}_k$ yields
\begin{equation}
p(\vec{Y}_k|\vec{\theta})=p(\vec{y}_k|\vec{Y}_{k-1}, \vec{\theta})p(\vec{Y}_{k-1}|\vec{\theta}), 
\label{eq:lem3p_Yk_recursion}
\end{equation}
where
\begin{equation}
p(\vec{y}_k|\vec{Y}_{k-1},\vec{\theta})=\int p(\vec{y}_k|\vec{x}_k, \vec{\theta})p(\vec{x}_k| \vec{Y}_{k-1}, \vec{\theta})d\vec{x}_k,
\label{eq:lem3p_yk_pred_distr}
\end{equation}
is the predictive distribution of $\vec{y}_k$. Dividing \eqref{eq:lem3p_filter_der} by \eqref{eq:lem3p_Yk_recursion}, gives the correction step:
\begin{equation}
p(\vec{x}_k|\vec{Y}_k,\vec{\theta})=\frac{p(\vec{y}_k|\vec{x}_k, \vec{\theta})p(\vec{x}_k| \vec{Y}_{k-1}, \vec{\theta})}{p(\vec{y}_k|\vec{Y}_{k-1},\vec{\theta})}.
\label{eq:lem3p_corr_step}
\end{equation}
Summarizing the above considerations, we have the following algorithm.\\
1)~Set the initial conditions:
\begin{equation}
p(\vec{x}_0|\vec{Y}_{-1}, \vec{\theta})=p(\vec{x}_0|\vec{\theta}), p(\vec{Y}_{-1}|\vec{\theta})=1.
\label{eq:lem3p_algorithm_ini}
\end{equation} 
2)~For $k=0, 1, \dots, N$, calculate:
\begin{align}
&p(\vec{y}_k|\vec{Y}_{k-1},\vec{\theta})=\int p(\vec{y}_k|\vec{x}_k, \vec{\theta})p(\vec{x}_k| \vec{Y}_{k-1}, \vec{\theta})d\vec{x}_k,\\
\label{eq:lem3p_Yk_recursion_alg}
&p(\vec{Y}_k|\vec{\theta})=p(\vec{y}_k|\vec{Y}_{k-1}, \vec{\theta})p(\vec{Y}_{k-1}|\vec{\theta}),\\
&p(\vec{x}_k|\vec{Y}_k,\vec{\theta})=\frac{p(\vec{y}_k|\vec{x}_k, \vec{\theta})p(\vec{x}_k| \vec{Y}_{k-1}, \vec{\theta})}{p(\vec{y}_k|\vec{Y}_{k-1},\vec{\theta})d\vec{x}_k}, \\
&p(\vec{x}_{k+1}| \vec{Y}_{k}, \vec{\theta})=\int p(\vec{x}_{k+1}| \vec{x}_{k}, \vec{\theta})p(\vec{x}_k|\vec{Y}_k, \vec{\theta})d\vec{x}_k.
\label{eq:lem3p_algorithm}
\end{align}
By substituting \eqref{eq:lem3p_transition_pdf}, \eqref{eq:lem3p_y_pdf}, \eqref{eq:lem3p_xini_pdf} into \eqref{eq:lem3p_algorithm_ini}-\eqref{eq:lem3p_algorithm}, and after somewhat tedious calculations, we obtain:
\begin{align}
&p(\vec{x}_k|\vec{Y}_k,\vec{\theta})=\cN(\vec{x}_k,\vec{m}_k(\vec{\theta}), \vec{S}_k(\vec{\theta})),\\
&p(\vec{y}_k|\vec{Y}_{k-1},\vec{\theta})=\cN(\vec{y}_k,\mat{C}\vec{m}_k^-(\vec{\theta}), \vec{\Sigma}_k(\vec{\theta})),
\label{eq:lem3p_yk_pred_pdf}
\end{align}
where $\vec{m}_k(\vec{\theta})$, $\vec{S}_k(\vec{\theta})$, $\vec{m}_k^-(\vec{\theta})$, $\vec{\Sigma}_k(\vec{\theta})$, $k=0,1,\dots$, are given by the Kalman filtering equations \eqref{eq:pem_kalman_0}-\eqref{eq:pem_kalman}. Then, by using the recursive formula \eqref{eq:lem3p_Yk_recursion_alg} and \eqref{eq:lem3p_yk_pred_pdf}, we get:
\begin{equation}
p(\vec{Y}|\vec{\theta})=\prod\limits_{k=0}^N\cN(\vec{y}_k,\mat{C}\vec{m}_k^-(\vec{\theta}), \vec{\Sigma}_k(\vec{\theta})),
\label{eq:lem3p_Y_pdf_formula_prod}
\end{equation}
where $\vec{Y}=\mrm{col}(\vec{y}_0, \dots, \vec{y}_N)$. On the other side, accordingly to \eqref{eq:model_lin_sys_0}-\eqref{eq:big_f_def_quasi_lin}, we have:
\begin{equation}
p(\vec{Y}|\vec{\theta})=\cN(\vec{Y},\vec{F}(\vec{\theta}), \mat{S}(\vec{\theta}))=\prod\limits_{k=0}^N\cN(\vec{y}_k,\mat{C}\vec{m}_k^-(\vec{\theta}), \vec{\Sigma}_k(\vec{\theta})).
\label{eq:lem3p_pdfs_compare}
\end{equation}
Taking the logarithm of both sides and calculating its expected value we get:
\begin{equation}
\begin{split}
&\int p(\vec{Y}|\vec{\theta})\left(-\ln p(\vec{Y}|\vec{\theta})\right)d\vec{Y}=\\
&=\tfrac{1}{2}\ln \left((2\pi e)^{n_y(N+1)} |\mat{S}(\vec{\theta})|\right)=\tfrac{1}{2}\ln \left((2\pi e)^{n_y(N+1)} \prod\limits_{k=0}^N|\vec{\Sigma}_k(\vec{\theta})|\right).
\end{split}
\label{eq:lem3p_log_py_exp}
\end{equation}
Hence $|\mat{S}|=\prod\limits_{k=0}^N|\vec{\Sigma}_k|$, which proves \eqref{eq:det_S_calc}. Now, taking the logarithm of \eqref{eq:lem3p_pdfs_compare}, we have:
\begin{equation}
\tfrac{1}{2}\ln|\mat{S}|+\tfrac{1}{2}|\vec{Y}-\vec{F}|_{\mat{S}^{-1}}^2=\frac{1}{2}\ln \left(\prod\limits_{k=0}^N|\vec{\Sigma}_k|\right)+\tfrac{1}{2}\sum\limits_{k=0}^N|\vec{y}_k-\mat{C}\vec{m}_k^-|_{\vec{\Sigma}_k^{-1}}^2,
\end{equation}
where the arguments has been omitted for convenience. Since $|\mat{S}|=\prod\limits_{k=0}^N|\vec{\Sigma}_k|$,
we get \eqref{eq:norm_Y_F_calc}. Putting \eqref{eq:norm_Y_F_calc} and \eqref{eq:det_S_calc} into \eqref{eq:log_likelihood_0} gives \eqref{eq:log_likelihood_recursive}.
\end{proof}
%%%% Lemma 4 proof ***************************************************************************************************
\begin{proof}[Proof of lemma 4]
Let $\vec{\Theta}=\lbrace\vec{\theta}_1,...,\vec{\theta}_r\rbrace, \vec{\theta}_i \in \mathbb{R}^{n_{\theta}}$. We are interested in calculation of the quantity 
\begin{equation}
d_{i,j}(\vec{U})=\tfrac{1}{8}\vec{\Delta}_{i,j}^{\TT}\left(\tfrac{1}{2}(\mat{S}_i+\mat{S}_j)\right)^{-1}\vec{\Delta}_{i,j}+\tfrac{1}{2}\ln |\tfrac{1}{2}(\mat{S}_i+\mat{S}_j)|-\tfrac{1}{4}\ln\left(|\mat{S}_i||\mat{S}_j|\right),
\label{eq:d_ij_proof_lem_4}
\end{equation}
where
\begin{equation}
    \vec{\Delta}_{i,j}=\vec{F}(\vec{\theta}_i, \vec{U})-\vec{F}(\vec{\theta}_j, \vec{U}), \mat{S}_i=\mat{S}(\vec{\theta}_i, \vec{U}), \mat{S}_j=\mat{S}(\vec{\theta}_j, \vec{U})
\end{equation}
and $\vec{F}(\vec{\theta}, \vec{U})$, $\mat{S}(\vec{\theta}, \vec{U})$, are defined by \eqnsref{eq:model_lin_sys_0}{eq:big_f_def_quasi_lin} of Sect. \ref{sec:quasi_linear_sys}. Let us define
\begin{align}
\label{eq:Y_ij_def_0}
&\vec{Y}^{(i)}=\vec{F}(\vec{\theta}_i, \vec{U})+\vec{Z}^{(i)}, \vec{Y}^{(j)}=\vec{F}(\vec{\theta}_j, \vec{U})+\vec{Z}^{(j)},\\
&\tilde{\vec{Y}}=\tfrac{1}{\sqrt{2}}\left(\vec{Y}^{(i)}-\vec{Y}^{(j)}\right)=\tfrac{1}{\sqrt{2}}\left(\vec{\Delta_{i,j}}+\vec{Z}^{(i)}-\vec{Z}^{(j)}\right),
\label{eq:Y_ij_def}
\end{align}
where $\vec{Z}^{(i)}\sim\cN(0,\mat{S}_i)$, $\vec{Z}^{(j)}\sim\cN(0,\mat{S}_j)$. The density of variable $\tilde{\vec{Y}}$, given $\tilde{\vec{\theta}}=\mrm{col}(\vec{\theta_i}, \vec{\theta_j})$, has the form
\begin{equation}
    p(\tilde{\vec{Y}}| \tilde{\vec{\theta}})=\cN\left(\tilde{\vec{Y}},\tfrac{1}{\sqrt{2}}\vec{\Delta_{i,j}},\tfrac{1}{2}(\mat{S}_i+\mat{S}_j)\right).
    \label{eq:big_y_tilde_pdf}
\end{equation}
Now, let us consider the following two systems:
\begin{align}
\label{eq:sys_2n_dynamics_0}
&\vec{x}^{(i)}_{k+1} =\mat{A}(\vec{\theta}_i, \vec{u}_k)\vec{x}^{(i)}_{k}+\mat{B}(\vec{\theta}_i, \vec{u}_k)+\mat{G} (\vec{\theta}_i, \vec{u}_k)\vec{w}^{(i)}_{k}, \vec{y}^{(i)}_{k} =\mat{C}\vec{x}^{(i)}_{k}+\vec{v}^{(i)}_k,\\
&\vec{x}^{(j)}_{k+1} =\mat{A}(\vec{\theta}_j, \vec{u}_k)\vec{x}^{(j)}_{k}+\mat{B}(\vec{\theta}_j, \vec{u}_k)+\mat{G}(\vec{\theta}_j, \vec{u}_k)\vec{w}^{(j)}_{k}, \vec{y}^{(j)}_{k} =\mat{C}\vec{x}^{(j)}_{k}+\vec{v}^{(j)}_k,
\label{eq:sys_2n_dynamics}
\end{align}
where $\vec{w}^{(i)}_{k}, \vec{w}^{(j)}_{k}\sim\cN(0,\mathbb{I})$ and $\vec{v}^{(i)}_{k}, \vec{v}^{(j)}_{k}\sim\cN(0,\mat{S}_v)$ are mutually independent. The components of the vectors $\vec{Y}^{(i)}$ and $\vec{Y}^{(i)}$ in \eqref{eq:Y_ij_def_0} correspond to the outputs of the systems \eqref{eq:sys_2n_dynamics_0} and \eqref{eq:sys_2n_dynamics}, that is $\vec{Y}^{(i)}=\mrm{col}(\vec{y}^{(i)}_0, \dots, \vec{y}^{(i)}_N)$, $\vec{Y}^{(j)}=\mrm{col}(\vec{y^{(j)}}_0, \dots, \vec{y}^{(j)}_N)$.
Then on the basis of \eqref{eq:Y_ij_def}, we get $\tilde{\vec{Y}}=\mrm{col}(\tilde{\vec{y}}_0, \dots, \tilde{\vec{y}}_0)$, where
\begin{equation}
    \tilde{\vec{y}}_k=\tfrac{1}{\sqrt{2}}\left(\vec{y}^{(i)}_k-\vec{y}^{(j)}_k\right)=\tfrac{1}{\sqrt{2}}\left(\mat{C}\left(\vec{x}^{(i)}_k-\vec{x}^{(j)}_k\right)+\left(\vec{v}^{(i)}_k-\vec{v}^{(j)}_k\right)\right).
    \label{eq:tilde_y}
\end{equation}
Defining the matrices $\tilde{\mat{A}}_k$, $\tilde{\mat{B}}_k$, $\tilde{\mat{G}}_k$, $\tilde{\mat{C}}$, according to \eqref{eq:system_2n_ab} and \eqref{eq:system_2n_gc}, and taking into account that $\tfrac{1}{\sqrt{2}}\left(\vec{v}^{(i)}_k-\vec{v}^{(j)}_k\right)\sim\cN(0, \mat{S}_v)$, we can replace \eqnsref{eq:sys_2n_dynamics_0}{eq:tilde_y} with a single, 
$2n$-dimensional system with $n_y$ outputs: 
\begin{equation}
    \tilde{\vec{x}}_{k+1}=\tilde{\mat{A}}_k(\tilde{\vec{\theta}})\tilde{\vec{x}}_k+\tilde{\mat{B}}_k(\tilde{\vec{\theta}})+\tilde{\mat{G}}_k(\tilde{\vec{\theta}})\tilde{\vec{w}}_k, \tilde{\vec{y}}_k=\tilde{\mat{C}}\tilde{\vec{x}}_k+\vec{v}_k,
\end{equation}
where $\tilde{\vec{x}}_k=\mrm{col}(\vec{x}^{(i)}_{k}, \vec{x}^{(j)}_{k})$, $\tilde{\vec{w}}_k=\mrm{col}(\vec{w}^{(i)}_{k}, \vec{w}^{(j)}_{k})$ and $\vec{v}_k\sim\cN(0,\mat{S}_v)$. Proceeding analogously to the proof of Lemma \ref{lem:Lemma_3}, we infer that the conditional density of variable $\tilde{\vec{Y}}$ is given by:
\begin{equation}
p(\tilde{\vec{Y}}|\tilde{\vec{\theta}})=\prod\limits_{k=0}^N\cN(\tilde{\vec{y}}_k, \tilde{\mat{C}}\tilde{\vec{m}}_k^-(\tilde{\vec{\theta}}), \tilde{\vec{\Sigma}}_k(\tilde{\vec{\theta}})),
\label{eq:small_y_tilde_pdf}
\end{equation}
where $\tilde{\vec{m}}_k^-(\tilde{\vec{\theta}})$, $\tilde{\vec{\Sigma}}_k(\tilde{\vec{\theta}})$, are calculated recursively by the Kalman filter equations
\begin{align}
\label{eq:2n_kalman_0}
&\tilde{\vec{\Sigma}}_k=\mat{S}_v+\tilde{\mat{C}}\tilde{\mat{S}}_k^-\tilde{\mat{C}}^{\TT},\\
&\tilde{\mat{L}}_k=\tilde{\mat{S}}_k^-\tilde{\mat{C}}^{\TT}\tilde{\vec{\Sigma}}_k^{-1},\\
\label{eq:2n_kalman_1}
&\tilde{\vec{m}}_k=(\mathbb{I}-\tilde{\mat{L}}_k\tilde{\mat{C}})\tilde{\vec{m}}_k^-+\tilde{\mat{L}}_k\tilde{\vec{y}}_k,\\
&\tilde{\mat{S}}_k=\tilde{\mat{S}}_k^--\tilde{\mat{L}}_k\tilde{\vec{\Sigma}}_k\tilde{\mat{L}}_k^{\TT},\\
&\tilde{\vec{m}}_{k+1}^-=\tilde{\mat{A}}_{k}\tilde{\vec{m}}_{k}+\tilde{\mat{B}}_{k},\\
&\tilde{\mat{S}}_{k+1}^-=\tilde{\mat{A}}_{k}\tilde{\mat{S}}_{k}\tilde{\mat{A}}_{k}^{\TT}+\tilde{\mat{G}}_{k}\tilde{\mat{G}}_{k}^{\TT}, k=0, 1 ..., N,
\label{eq:2n_kalman}
\end{align}
with initial conditions \eqref{eq:init_cond_d_ij}. Comparing \eqref{eq:big_y_tilde_pdf} with \eqref{eq:small_y_tilde_pdf} gives:
\begin{equation}
    p(\tilde{\vec{Y}}| \tilde{\vec{\theta}})=\cN\left(\tilde{\vec{Y}},\tfrac{1}{\sqrt{2}}\Delta_{i,j},\tfrac{1}{2}(\mat{S}_i+\mat{S}_j)\right)=\prod\limits_{k=0}^N\cN(\tilde{\vec{y}}_k, \tilde{\mat{C}}\tilde{\vec{m}}_k^-(\tilde{\vec{\theta}}), \tilde{\vec{\Sigma}}_k(\tilde{\vec{\theta}})).
    \label{eq:y_tilde_pdf_compare}
\end{equation}
Taking the logarithm of both sides and calculating its expected value we get:
\begin{equation}
\begin{split}
&\int p(\tilde{\vec{Y}}|\tilde{\vec{\theta}})\left(-\ln p(\tilde{\vec{Y}}|\tilde{\vec{\theta}})\right)d\tilde{\vec{Y}}=\\
&=\tfrac{1}{2}\ln \left((2\pi e)^{n_y(N+1)} |\tfrac{1}{2}(\mat{S}_i+\mat{S}_j)|\right)=\tfrac{1}{2}\ln \left((2\pi e)^{n_y(N+1)} \prod\limits_{k=0}^N|\tilde{\vec{\Sigma}}_k(\vec{\tilde{\theta}})|\right).
\end{split}
\label{eq:lem4p_log_py_exp}
\end{equation}
Hence
\begin{equation}
\ln|\tfrac{1}{2}(\mat{S}_i+\mat{S}_j)|=\sum\limits_{k=0}^N\ln|\tilde{\vec{\Sigma}}_k|.
\label{eq:lem4p_dets_equal}
\end{equation}
Taking the logarithm of \eqref{eq:y_tilde_pdf_compare} and applying \eqref{eq:lem4p_dets_equal} yields:
\begin{equation}
\tfrac{1}{2}(\tilde{\vec{Y}}-\tfrac{1}{\sqrt{2}}\vec{\Delta_{i,j}})^{\TT}\left(\tfrac{1}{2}(\mat{S}_i+\mat{S}_j)\right)^{-1}(\tilde{\vec{Y}}-\tfrac{1}{\sqrt{2}}\vec{\Delta_{i,j}})=\tfrac{1}{2}\sum\limits_{k=0}^N|\tilde{\vec{y}}_k-\tilde{\mat{C}}\tilde{\vec{m}}_k^-|_{\tilde{\vec{\Sigma}}_k^{-1}}^2.
\label{eq:lem4p_norms_equal}
\end{equation}
Finally, since $\tilde{\vec{Y}}=\mrm{col}(\tilde{\vec{y}}_0, \dots, \tilde{\vec{y}}_0)$, then substituting $\tilde{\vec{Y}}=0$, $\tilde{\vec{y}}_k=0$ in \eqref{eq:lem4p_norms_equal}, \eqref{eq:2n_kalman_1}, we conclude that:  
\begin{equation}
\tfrac{1}{8}\Delta_{i,j}^{\TT}\left(\tfrac{1}{2}(\mat{S}_i+\mat{S}_j)\right)^{-1}\Delta_{i,j}=\tfrac{1}{4}\sum\limits_{k=0}^N|\tilde{\mat{C}}\tilde{\vec{m}}_k^-|_{\tilde{\vec{\Sigma}}_k^{-1}}^2.
\end{equation}
where $\tilde{\vec{m}}_k^-$, $\tilde{\vec{\Sigma}}_k$ fulfils equations \eqref{eq:kalman_ij_0}-\eqref{eq:kalman_ij_1}. The last term in \eqref{eq:d_ij_proof_lem_4} is calculated according to Lemma \ref{lem:Lemma_3}. 
\end{proof}
%*****************************************************
\section{An example of the gap between ITB and BCRB}
\label{app:appendix_B}
The difference between ITB \eqref{eq:error_lb} and BCRB \eqref{eq:van_trees} can be significant. To see this, let us consider the model
\begin{equation}
    y=\theta+v,
    \label{eq:efrom_example_model}
\end{equation}
where $v\sim\cN(0, 1)$. The elementary calculation yields $J_D=1$. If prior is Gaussian, that is, $p_0(\theta)=\cN(\theta, m_{\theta}, \sigma_{\theta}^2)$, then $J_P=\sigma_{\theta}^{-2}$, $2(H_{\theta}-I_{\theta;y})=\ln 2\pi e-\ln(\sigma_{\theta}^{-2}+1)$ and both BCRB \eqref{eq:van_trees} and ITB \eqref{eq:error_lb} yield the same error estimate. Now, let us assume that
\begin{equation}
    p_0(\theta)=\tfrac{1}{2}\left(\Phi(\alpha(1+\theta))+\Phi(\alpha(1-\theta))-1\right),
    \label{eq:efrom_example_prior}
\end{equation}
where $\alpha>0$ is a parameter and $\Phi(t)=\int\limits_{-\infty}^{t}\cN (t,0,1)dt$. The prior \eqref{eq:efrom_example_prior} is an analytic function which, in the limit $\alpha\rightarrow\infty$, tends to the uniform distribution $U[-1,1]$. Since $H(y|\theta)=\tfrac{1}{2}\ln(2\pi e)$, $H(y)\leqslant\tfrac{1}{2}\ln(2\pi e(\mrm{var}(\theta)+1))$, $\mrm{var}(\theta)=\tfrac{1}{3}+O_1(\alpha^{-1})$, $H(\theta)=\ln2+O_2(\alpha^{-1})$ and $H_{\theta}-I_{\theta;y}=H(y|\theta)+H(\theta)-H(y)$, after elementary calculations, we get the ITB:
\begin{equation}
    \E(\theta-\hat{\theta}(y))^2\geqslant\frac{e^{2(H_{\theta}-I_{\theta;y})}}{2\pi e}\geqslant\frac{3\sqrt{2}}{8\pi e}+O(\alpha^{-1}).
\end{equation}
On the other hand, according to \eqref{eq:J_P}, we have:
\begin{equation}
    J_P=\frac{\alpha^2}{2}\int\limits_{-\infty}^{\infty}\frac{\left(\cN(\alpha(1+\theta),0,1)-\cN(\alpha(1-\theta),0,1)\right)^2}{\Phi(\alpha(1+\theta))+\Phi(\alpha(1-\theta))-1}d\theta\geqslant\frac{\alpha}{2\sqrt{\pi}}  \xrightarrow[\alpha\rightarrow\infty]{}\infty.
\end{equation}
Hence, BCRB \eqref{eq:van_trees} becomes trivial, but ITB still gives a reasonable error estimate. If the likelihood is non-Gaussian, a similar effect occurs. Therefore, BCRB generally underestimates the minimum possible estimation error.
\section{Discretization of linear SDE}
\label{app:appendix_C}
Consider continuous-time SDE
\begin{equation}
d\vec{x}=(\mat{A}_c\vec{x}+\mat{B}_c\vec{u})dt+\mat{G}_cd\vec{w},
\label{eq:ct_sde_app_b}
\end{equation}
where $\vec{x}(t)\in\mathbb{R}^n$, $\vec{w}(t)\in\mathbb{R}^{n_w}$ is a vector of mutually independent standard Wiener processes. Let $\Delta$ denote the discretization period, and let $\vec{u}(t)=\vec{u}_k$, $t\in [t_k, t_{k+1}]$, $t_k=k\Delta$. Then process $\vec{x}_k=\vec{x}(t_k)$, fulfils the difference equation 
\begin{equation}
\vec{x}_{k+1}=\mat{A}\vec{x}_k+\mat{B}u_k+\mat{G}\vec{w}_k,
\end{equation}
where $\vec{w}_k\sim\cN(0,\mathbb{I}_{n_w})$ and 
\begin{equation}
\mat{A}=e^{\mat{A}_{c}\Delta}, \mat{B}=\int\limits_0^{\Delta}e^{\mat{A}_{c}\tau}\mat{B}_cd\tau, \mat{D}=\mat{G}\mat{G}^{\TT}=\int\limits_0^{\Delta}e^{\mat{A}_{c}\tau}\mat{G}_c\mat{G}_c^{\TT}e^{\mat{A}_{c}^{\TT}\tau}d\tau. 
\end{equation}
If $\mat{D}>0$, then $\mat{G}$ can be determined by the Cholesky factorization of $\mat{D}$. In the general case we use spectral decomposition $\mat{D}=\mat{Q}\vec{\Lambda}\mat{Q}^{\TT}$, and then $\mat{G}=\mat{Q}\vec{\Lambda}^{0.5}$.
\end{appendices}
\vspace{0.5cm}

\textbf{Abbreviations.}
The following abbreviations are used in this manuscript: FIM - Fisher Information Matrix; CRB - Cram\'{e}r-Rao
 Bound; BCRB - Bayesian Cram\'{e}r-Rao
 Bound; ITB - Information-Theoretic Bound; DOE - Design of Experiment; SDE - Stochastic Differential Equation. The norm of the vector $x\in\mathbb{R}^n$ is denoted by $|x|$. For any square matrix $\mat{Q}$, the quadratic form $\vec{x}^{\TT}\mat{Q}\vec{x}$ is denoted by $|\vec{x}|_\mat{Q}^2$. The trace and determinant of the matrix $\mat{A}$ are denoted by $\mrm{tr}\mat{A}$, $|\mat{A}|$ or $\det(\mat{A})$. The set of symmetric, positive definite matrices of dimension \textit{n} is denoted by $\mat{S}^{+}(n)$. The symbol $\mathrm{col}(\vec{a}_{1} ,\vec{a}_{2} ,...,\vec{a}_{n})$ denotes the column vector. $\vec{\xi} \sim \cN(\vec{m},\mat{S}) $ means that $\vec{\xi}$ has a normal distribution with mean $\vec{m}$ and covariance $\mat{S}$. The density of the Gaussian variable is denoted by $\cN(\vec{x},\vec{m},\mat{S})=2\pi ^{ -\frac{n}{2}} |\mat{S}|^{ -\frac{1}{2}} \exp (-0.5(\vec{x}-\vec{m})^{\TT} \mat{S}^{ -1}(\vec{x}-\vec{m}))$.
\vspace{0.5cm}

\textbf{Author contributions.}
Article concept, proofs of Theorem \ref{thm:theorem_1} and Lemmas \ref{lem:Lemma_2}, \ref{lem:Lemma_3}, \ref{lem:lemma_4}, all Appendices, MATLAB codes, implementation of Bayesian and classical signal selection methods, development of all examples, comparison with classical methods, all formulas derivations and text writing: Piotr Bania. Determination of the optimal spectrum and signal generation using the MOOSE-2 toolbox in Sect. \ref{sec:example_01}, \ref{sec:example_02}, verification of the correctness of formulas describing discrete systems in Sect. \ref{sec:example_01}, \ref{sec:example_02}, \ref{sec:example_03} and verification of the proofs of Lemmas 3 and 4: Anna W{\'o}jcik. Text proofreading, literature review, introduction, discussion and conclusions: Piotr Bania and Anna W{\'o}jcik.
\vspace{0.5cm}

\textbf{Funding.} This research was supported by the statutory subsidy of the AGH University of Science and Technology,
No. 16.16.120.773 and by the Initiative of Excellence - Research University (IDUB) program.
\vspace{0.5cm}

\textbf{Data availability.} The MATLAB codes, in particular the functions for calculating the lower bound \eqref{eq:inf_lb_formulas_0} and the $d_{i,j}$ in \eqref{eq:inf_lb_formulas_dij}, \eqref{eq:opt_problem_two_poss}, are available on \href{https://github.com/Jhiqo/Bay_design_ql_sys}{GitHub} repository.
\vspace{0.5cm}

\textbf{Acknowledgments.} The authors gratefully acknowledge Morgan Mitchell, Jan Ko{\l}ody{\'n}ski, Klaudia Dilcher, Aleksandra Sierant, Julia Amor{\'o}s-Binefa, and Diana M{\'e}ndez-{\'A}valos for insightful discussions on quantum control and atomic magnetometers.
\end{document}